\documentclass[journal]{IEEEtran}

\usepackage[T1]{fontenc}
\usepackage{graphicx}
\usepackage[caption=false]{subfig}
\usepackage{mathrsfs,amsmath,amsfonts,amssymb,amsthm,amstext,amscd,amsxtra,amsopn}
\usepackage{diagbox}
\usepackage{lettrine}
\usepackage{cite}
\usepackage{eucal}
\providecommand{\keywords}[1]{\textbf{\textit{Index terms---}} #1}
\usepackage{bm}
\usepackage{multirow}
\usepackage{esint}
\usepackage{empheq}
\usepackage{bm}
\usepackage{algorithm}
\usepackage[noend]{algpseudocode}
\usepackage{etoolbox}
\usepackage{leftidx}
\usepackage{makecell}
\usepackage{stackrel}
\usepackage{xcolor}
\usepackage{nicefrac}
\usepackage[normalem]{ulem}
\usepackage[hyphens]{url}
\usepackage[hidelinks]{hyperref}
\hypersetup{breaklinks=true}
\renewcommand{\vec}[1]{\bm{#1}}

\newcommand{\matvec}[1]{\mathbf{#1}}

\def\J{\vec{j}}
\def\M{\vec{m}}
\def\E{\vec{e}}
\def\H{\vec{h}}

\def\Ei{\vec{e}_{\rm inc}}
\def\Hi{\vec{h}_{\rm inc}}

\DeclarePairedDelimiter\abs{\lvert}{\rvert}
\newcommand\norm[1]{\left\lVert#1\right\rVert}

\def\ce{{c_e}}
\def\cm{{c_m}}

\newcommand\Aop{{\CMcal{A}}}
\newcommand\Bop{{\CMcal{B}}}
\newcommand\Eop{{\CMcal{E}}}

\newcommand\Gop{{\CMcal{G}}}

\newcommand\Kop{{\CMcal{K}}}

\newcommand\Nop{{\CMcal{N}}}
\newcommand\Oop{{\CMcal{O}}}

\newcommand\At{{\cal{A}}}

\newcommand\Kt{{\cal{K}}}

\newcommand\Mt{{\cal{M}}}
\newcommand\Nt{{\cal{N}}}

\newcommand\chiE{{\chi_e}}
\newcommand\chiH{{\chi_m}}

\makeatletter
\def\BState{\State\hskip-\ALG@thistlm}

\expandafter\patchcmd\csname\string\algorithmic\endcsname{\itemsep\z@}{\itemsep=1.5mm plus2pt}{}{}

\newenvironment{breakablealgorithm}
  {
   \begin{center}
     \refstepcounter{algorithm}
     \hrule height.8pt depth0pt \kern2pt
     \renewcommand{\caption}[2][\relax]{
       {\raggedright\textbf{\ALG@name~\thealgorithm} ##2\par}%
       \ifx\relax##1\relax 
         \addcontentsline{loa}{algorithm}{\protect\numberline{\thealgorithm}##2}%
       \else 
         \addcontentsline{loa}{algorithm}{\protect\numberline{\thealgorithm}##1}%
       \fi
       \kern2pt\hrule\kern2pt
     }
  }{
     \kern2pt\hrule\relax
   \end{center}
  }

\makeatother

\algnewcommand{\Initialize}[1]{%
  \State \textbf{Initialize:}
  \Statex \hspace*{\algorithmicindent}\parbox[t]{.8\linewidth}{\raggedright #1}
}

\newcommand{\iu}{{\mathrm{i}}}

\title{Memory footprint reduction for the \\ FFT-based volume integral equation method \\ via tensor decompositions}

\author{Ilias I. Giannakopoulos, Mikhail S. Litsarev and Athanasios G. Polimeridis ~\IEEEmembership{Senior Member,~IEEE}

\thanks{Ilias I. Giannakopoulos and Mikhail S. Litsarev are with the Skoltech Center for Computational Data-Intensive Science and Engineering, Skolkovo Institute of Science and Technology, 143026 Moscow, Russia.}
\thanks{Athanasios G. Polimeridis is with Q Bio, CA 94063, USA.}}

\begin{document}
\bstctlcite{IEEEexample:BSTcontrol}

\maketitle

\begin{abstract}
We present a method of memory footprint reduction for FFT-based, electromagnetic (EM) volume integral equation (VIE) formulations. The arising Green's function tensors have low multilinear rank, which allows Tucker decomposition to be employed for their compression, thereby greatly reducing the required memory storage for numerical simulations. Consequently, the compressed components are able to fit inside a graphical processing unit (GPU) on which highly parallelized computations can vastly accelerate the iterative solution of the arising linear system. In addition, the element-wise products throughout the iterative solver's process require additional flops, thus, we provide a variety of novel and efficient methods that maintain the linear complexity of the classic element-wise product with an additional multiplicative small constant. We demonstrate the utility of our approach via its application to VIE simulations for the Magnetic Resonance Imaging (MRI) of a human head. For these simulations we report an order of magnitude acceleration over standard techniques.    
\end{abstract}


\section{Introduction} \label{sc:I}

Magnetic resonance imaging has admittedly become an essential tool of modern medical imaging and disease diagnosis. Indeed, modern MR scanners offer extremely precise and detailed views of the human body. Nevertheless, the MR community is constantly looking for new ideas that will improve further their capabilities, especially focusing on higher resolution and shorter scanning times, both directly related to the signal-to-noise ratio (SNR). The current trend for increasing the SNR of the MR scanners is moving to higher magnetic fields; only in the last decade the advances in superconducting magnet technology allowed to move from 1.5 to 7 Tesla clinical scanners \cite{siemens2017seventesla}. However, as the strength of the magnetic field increases, so does the operating radio-frequency (RF) of the associated coils. Hence the interactions between electromagnetic waves and biological tissue become more dominant, and can easily have detrimental effect to the quality of the images and the safety of the patients, if not modeled accurately while designing the scanner and the RF coils \cite{Lattanzi2009, jin1997sar, zhang2014predicting, cosottini2014short}.
\par
A plethora of well-established methods of computational electromagnetics has been extensively used over the past decades to study the time-harmonic solutions of Maxwell's equations. On the one hand, partially differential equation methods, such as the finite difference and the finite element method are a great tool for the EM simulation of inhomogeneous arbitrary shaped objects \cite{dimbylow1997fdtd, collins2001calculations, kuhn2009assessment}. On the other hand, when the maximal use of a specific setting is possible, the VIE method provides an opportunity for customization of fast algorithms. Towards that direction, a magnetic-resonance integral equation-based suite (MARIE) \cite{MARIE, Polimeridis2014, Villena2016} has been developed, where polynomial basis functions are used for the fast and precise EM modeling of the interactions between human tissue and MR coils. Specifically, by expanding the unknowns with higher order polynomials \cite{georgakis2019fast}, superior numerical accuracy is in place, contrary to standard low-order approximations, even for the challenging dielectric shimming technique \cite{brink2014high, koolstra2018improved}.
\par
When VIEs are discretized over a uniform grid with polynomial basis functions, the arising Galerkin Method of Moments (MoM) system matrix has block-Toeplitz with Toeplitz-blocks (BTTB) structure, owing to the translation invariance property of the Green's function. Therefore, only the defining columns (formulated as tensors) need to be stored for the matrix-vector product implementation, which can be greatly accelerated with the help of the FFT algorithm \cite{Borup1984, Schaubert1984, Catedra1989, Zwamborn1992, Zwamborn1994, Gan1994, Jin1996, Ozdemir2006, Ozdemir2007, Beurden2008, Markkanen2012, Markkanen2012b, Oijala2014}. Regrettably, when we discretize FFT-based VIEs with higher order polynomials, namely piecewise linear (PWL), the required memory footprint of the arising Green's function tensors increases significantly, thus, forbidding the usage of heterogeneous computing techniques; clearly, for large computational domains, the aforementioned tensors might be prohibitively large for storing in fast memory. However, for MRI applications, where the dimensions of the computational domain (human body) are comparable with the operating wavelength, the Green's function volume integral operators present low rank properties \cite{Chai2013}. Hence, low multilinear rank tensor decompositions can be applied to dramatically reduce the required storage memory. 
\par
Tensor decompositions are commonly used tools for the analysis of multidimensional data in various scientific fields \cite{Cichocki2009, Sidiropoulos2000b, Lathauwer2007, Comon2014}. Specifically, Tucker decomposition provides an optimal fit and a stable approximation for three-dimensional (3D) tensors \cite{Tucker1966}. It can be implemented either with the well-known higher-order singular value decomposition (HOSVD) \cite{Lathauwer2000}, which has been successfully applied in multidimensional data analysis in the past years \cite{Vasilescu2002, Duchenne2011, Rajwade2013, Yucel2017}, or with cross approximation-based techniques \cite{Oseledets2008, Rakhuba2015}. Another interesting decomposition is the canonical polyadic (CP) model \cite{Hitchcock1927a, Hitchcock1927b, Carroll1970, Harshman1970}, which gives the most compact representation of the initial array. On contrast with Tucker algorithms, CP algorithms are often ill-posed for big tensors, meaning that a low multilinear rank approximation might not exist at all \cite{Silva2008}. To avoid the aforementioned impasse the Tucker+CP decomposition \cite{Bro1998} can be utilized, where the CP is applied on a small full tensor (Tucker core).
\par
Our purpose in this article is to expand on the work presented in \cite{Polimeridis2014b,Giannakopoulos2018}, where Tucker decomposition was used to compress the Green's function tensors of the FFT-based, current VIE formulation (FFT-JMVIE), expanded with polynomial basis functions. Here, we compress the aforementioned tensors of JMVIE, expanded with PWL functions and perform the numerical simulations via a novel heterogeneous computing technique: First, we create all the required tensors in CPU and compress them via Tucker decomposition algorithms. Secondly, the resulting compressed forms are transferred to GPU, which has much less storage memory than CPU and cannot fit all the uncompressed components along with the vector of the unknowns. The element-wise products that appear in the matrix-vector product implementation require additional flops since tensor decompressions are needed, thus, we propose novel methods that maintain the linear complexity of the element-wise product with an additional multiplicative small constant. The proposed matrix-vector product consists of FFTs and matrix multiplications, which can be greatly accelerated via the highly parallel architecture of GPU. It is worth noting that the methods presented herein are applicable for all FFT-based VIE formulations. 
\par
To validate the algorithms presented in this work, we provide the compression factor of the arising Green's function tensors of FFT-JMVIE, along with time measurements for the novel matrix-vector product implementations both in C\texttt{++} and CUDA. In addition, we present the relative error of the absorbed power and the RF transmit field in MR measurements, between the standard approach and the implementations with HOSVD and Tucker+CP decompositions, for the case of a highly inhomogeneous realistic human head model. The results are in excellent agreement and the GPU-accelerated solvers are an order of magnitude faster. Finally, a refinement analysis shows the scaling properties of the novel methods. The remainder of this paper is organized as follows: In Section \ref{sc:II} we briefly describe the JMVIE formulation and we set up the associated linear system with polynomial basis functions. In addition, a theoretical study on the low multilinear rank property of the arising $\Nop$ and $\Kop$ tensors (Green's function tensors) is provided. In Section \ref{sc:III} we present an overview of Tucker and canonical polyadic models and in Section \ref{sc:IV} of HOSVD along with the Tucker+CP algorithm. In Section \ref{sc:V} we analyze the computational complexity of the associated operations related to the matrix-vector product implementation when tensor decompositions are used. In Section \ref{sc:VI}, the aforementioned results are presented. Finally, Table I lists some frequently used notation in this work.

\begin{table}[ht]
\caption{Notation} \label{tb:n1} \centering
\begin{tabular}{c| l }
\hline\hline\\[-0.4em]
Notation                 	  				 & Description                                                          \\[0.4em] \hline \\[-0.4em]
$a$                      	  				 & Scalar                                 	                            \\[0.3em]
$\vec{a}$                	  				 & Vector in $\mathbb{C}^3$                 		                    \\[0.3em]
$\At$              	          				 & Operator acting on vectors in $\mathbb{C}^3$                         \\[0.3em]
$\matvec{a}$              	  				 & Vector in $\mathbb{C}^{n}$                                           \\[0.3em]
$A$                      	  				 & Matrix in $\mathbb{C}^{n_1 \times n_2}$                              \\[0.3em]
$A^*$                    	  				 & Conjugate transpose of matrix                     				    \\[0.3em]
$A^{\left(\cdot \right)}$          			 & Two-dimensional reshape of $\Aop$                    				\\[0.3em]
$\Aop$                        			  	 & Tensor in $\mathbb{C}^{n_1 \times n_2 \times n_3}$                   \\[0.3em]
$\tilde{\cdot}$            				     & Approximation of array                                               \\[0.3em]
$\cdot_{i_1 i_2 i_3 \dots}$   				 & Element of array   	                                 		        \\[0.3em]
$\mathscr{F}\{\cdot\}$        			     & Discrete Fourier transform                                           \\[0.3em]
$\langle \cdot \rangle$, $\odot$, $\times_i$ & Inner, outer and n-mode products                                     \\[0.3em]
$\iu$                                        & Imaginary unit $\iu^2 = -1$                                          \\[0.3em]
\hline\hline
\end{tabular}
\end{table}

\section{Volume Integral Equations} \label{sc:II}

\subsection{Formulation}

We consider the 3D electromagnetic scattering by a closed penetrable object $\Omega\subset\mathbb{R}^3$. The permittivity and permeability of the scatterer are functions of $\vec{r}$, $\epsilon = \epsilon_r\left(\vec{r}\right)\epsilon_0$ and $\mu = \mu_r\left(\vec{r}\right)\mu_0$ respectively, where $\epsilon_0$ and $\mu_0$ are the vacuum electromagnetic properties, and $\epsilon_r$, $\mu_r$ are the relative properties of the scatterer. The wavelength is $\lambda$ and the wavenumber is $k_0 = \nicefrac{2 \pi}{\lambda}$. The JMVIE has the following form
\begin{equation}
\At \begin{bmatrix}
\J \\
\M 
\end{bmatrix} = \begin{bmatrix}
\ce & 0\\
0   & \cm
\end{bmatrix} \begin{bmatrix}
\Mt_{\chiE} & 0\\
0           & \Mt_{\chiH}
\end{bmatrix} \begin{bmatrix}
\Ei \\
\Hi 
\end{bmatrix}
\label{eq:n1}
\end{equation}
where
\begin{equation}
\begin{aligned}
\ce   \triangleq \iu \omega \epsilon_0 ,& \quad \cm   \triangleq \iu \omega \mu_0 \\
\chiE \triangleq \epsilon_r - 1      ,& \quad \chiH \triangleq \mu_r - 1.
\end{aligned}
\label{eq:n2}
\end{equation}
The electric and magnetic fields are denoted as $ \E,\H $ and the equivalent electric and magnetic polarization currents as
\begin{equation}
\begin{aligned}
\J &= \ce \left(\epsilon_r\left(\vec{r}\right)-1\right)\E\left(\vec{r}\right) \\
\M &= \cm \left(\mu_r\left(\vec{r}\right)-1\right)\H\left(\vec{r}\right).
\end{aligned}
\label{eq:n3}
\end{equation}
The operator $\At$ is 
\begin{equation}
\At = \begin{bmatrix}
\Mt_{\epsilon_r} - \Mt_{\chiE}\Nt & \ce\Mt_{\chiE}\Kt \\
-\cm\Mt_{\chiH}\Kt                & \Mt_{\mu_r} - \Mt_{\chiH}\Nt 
\end{bmatrix},
\label{eq:n4}
\end{equation}
the operators $\Mt_{\phi}$ denote a multiplication over the respective parameter $\phi$ and
\begin{equation}
\begin{aligned}
\Nt \left(\vec{s}\right) & \triangleq \nabla \times \nabla \times \int\limits_{\Omega} g\left(\vec{r}-\vec{r}'\right) \vec{s}\left(\vec{r}'\right)d^3\vec{r}' \\
\Kt \left(\vec{s}\right) & \triangleq \nabla \times \int\limits_{\Omega} g\left(\vec{r}-\vec{r}'\right) \vec{s}\left(\vec{r}'\right)d^3\vec{r}' 
\end{aligned}
\label{eq:n5}
\end{equation}
are the associated volume integral operators. The free-space Green's function, or fundamental Helmholtz solution, is
\begin{equation}
g\left(\vec{r} - \vec{r}' \right) = \frac{e^{-\iu k_0 \abs{\vec{r}-\vec{r}'}}}{4 \pi \abs{\vec{r}-\vec{r}'}}
\label{eq:n6}
\end{equation}
with $\vec{r}$ being the observation and $\vec{r}'$ the source point.

\subsection{Linear System}

In this work, we use the Galerkin MoM to numerically solve the aforementioned system of integral equation. Moreover, if the method is applied on uniform grids, with the number of voxels being $N_v = n_1 \times n_2 \times n_3$ ($n_1,n_2,n_3$ denote the discretization's linear size in the $x,y,z$ direction respectively), the arising matrices $N, K \in \mathbb{C}^{N_v \times N_v}$ admit a block-Toeplitz with Toeplitz blocks (BTTB) structure. We only need to store the Toeplitz defining tensors (Green's tensors) $\Nop \: \& \: \Kop$ $\in \mathbb{C}^{n_1 \times n_2 \times n_3}$ and then accelerate the associated matrix by vector products with the aid of the FFT algorithm along with an iterative solver. 
\par
The unknown volumetric currents should be expanded in some discrete set of appropriate basis functions e.g., PWC \cite{Polimeridis2014} or PWL \cite{georgakis2019fast} as follows:
\begin{equation}
j_q = \sum\limits_{m=1}^{N_v} \sum\limits_{l=1}^{L} u_m^{ql} C_m^{l}(\vec{r}).
\label{eq:n7a}
\end{equation}
The subscript $q \in \{x,y,z\}$ indicates the component of the current, $m$ corresponds to a specific voxel, geometrically defined by three indices $\{m_1,m_2,m_3\}$, with $(x_m,y_m,z_m)$ being its center, index $L$ is equal to $1$ for PWC and $4$ for PWL, $u_m^{ql}$ is the unknown current coefficient at the $m$-th voxel and $C_m^{q}$ are the scalar basis functions per voxel: 
\begin{equation}
\begin{aligned}
C_m^{1}(\vec{r}) &= P_m(\vec{r})\\
C_m^{2}(\vec{r}) &= \frac{x-x_m}{\Delta x}P_m(\vec{r})\\
C_m^{3}(\vec{r}) &= \frac{y-y_m}{\Delta y}P_m(\vec{r})\\
C_m^{4}(\vec{r}) &= \frac{z-z_m}{\Delta z}P_m(\vec{r})\\
\end{aligned}
\label{eq:n7b}
\end{equation}
where $P_m$ is a volumetric pulse, equal to $1$ inside the $m$-th voxel and $0$ elsewhere, and $\Delta x,\Delta y,\Delta z$ are the dimensions of each voxel in the $x,y,z$ directions respectively. In the case of PWC we are dealing with $3$ unknowns per voxel and in the case of PWL with $12$, therefore there are $9$ components per voxel each for the $\Nop$ and $\Kop$ tensors in the PWC approximation, and $144$ components per voxel for the $\Nop$ and $\Kop$ tensors in the PWL approximation. However, it can be seen for PWL that, due to symmetries and zero entries, there are only $60$ unique entries per voxel for the $\Nop$ tensor which need to be stored, and $30$ per voxel for the $\Kop$ tensor. Their elements are given via the standard $L_2$ inner product (defined in Appendix I) as follows 
\begin{equation}
\begin{aligned}
\Nop_{i j k}^{ql,q'l'} &= \langle \vec{f}_{1 1 1}^{ql}, \Nt \vec{f}_{i j k}^{q'l'} \rangle \\
\Kop_{i j k}^{ql,q'l'} &= \langle \vec{f}_{1 1 1}^{ql}, \Kt \vec{f}_{i j k}^{q'l'} \rangle
\end{aligned}
\label{eq:n8a}
\end{equation}
where
\begin{equation}
\begin{aligned}
\vec{f}_{ijk}^{q'l'} &= u_{\{i,j,k\}}^{q'l'} C_{\{i,j,k\}}^{l'}(\vec{r}) \hat{q'} \\
\vec{f}_{111}^{ql} &= u_{\{1,1,1\}}^{ql} C_{\{1,1,1\}}^{l}(\vec{r}) \hat{q}
\end{aligned}
\label{eq:n8b}
\end{equation}
are appropriate vector components of basis and testing functions in the voxel $\{i,j,k\}$ and $\{1,1,1\}$ respectively, and $\hat{q}, \hat{q}' \in \{\hat{x},\hat{y},\hat{z}\}$.    

\subsection{Ranks of \texorpdfstring{$\Nop$}{N} and \texorpdfstring{$\Kop$}{K} Tensors}

In \cite{Chai2013} it is proved that the Green's function related integrodifferential operators that arise from two well-separated geometry blocks have low-rank properties. Therefore, the off-diagonal blocks of the MoM matrix are low-rank since they represent such interactions. In the case of uniform grids, where only the defining BTTB tensor is stored, the interactions between one voxel and all the others are modeled as a column of the MoM matrix, meaning that we calculate interactions between distant voxels as well. Thus, a low multilinear rank approximation can be achieved in the tensor form of $\Nt$ and $\Kt$. According to \cite{Liu2012}, for 3D geometries, like the ones under study, the ranks of the operators, corresponding to interactions between well-separated geometrical source and observation domains, are proportional to the operating frequency $\mathcal{O}\left( k_0 \right)$. The case of study presented herein shares similarities in the dependence of the multilinear rank with frequency. However, since the voxel $\{1,1,1\}$ is part of both domains, due to the form of the Green's function tensors (\ref{eq:n8a}), not an exact linear dependence is expected. 

\section{Overview of Basic Tensor Decompositions} \label{sc:III}

A 3D array $\Aop \in \mathbb{C}^{n_1 \times n_2 \times n_3}$ can be approximated with a prescribed accuracy $\epsilon$ according to Tucker's model \cite{Tucker1966} as follows:
\begin{equation}
\Aop \triangleq \tilde{\Aop} + \Eop, \quad \norm{\Eop}_F = \epsilon
\label{eq:n12}
\end{equation}
where $\norm{\cdot}_F$ is the \textit{tensor Frobenius norm} (defined in Appendix I) and the approximation $\tilde{\Aop}$ is given from the following expression:
\begin{equation}
\tilde{\Aop}_{ijk} = \sum_{\alpha=1}^{r_1} \sum_{\beta=1}^{r_2} \sum_{\gamma=1}^{r_3} \Gop_{\alpha \beta \gamma} U^{1}_{i \alpha} U^{2}_{j \beta} U^{3}_{k \gamma}.
\label{eq:n13}
\end{equation}
The matrices $U^{q} \in \mathbb{C}^{n_q \times r_q}, q\text{=}1,2,3$ are the so-called \textit{Tucker factors} and the tensor $\Gop  \in \mathbb{C}^{r_1 \times r_2 \times r_3}$ is the \textit{Tucker core}. We can use the \textit{n-mode} products (also defined in Appendix I) and get the following compact form:
\begin{equation}
\tilde{\Aop} = \Gop \times_1 U^{1} \times_2 U^{2} \times_3 U^{3}.
\label{eq:n14}
\end{equation}
The triplet $\{r_1,r_2,r_3\}$ is called \textit{Tucker rank} or \textit{multilinear rank} of $\Aop$. A graphical representation of the Tucker decomposition of $\Aop$ is visualized in Fig. 1.(a). 
\par
For the compression of low-multilinear rank $d$-dimensional arrays with $d > 3$, the tensor train (TT) decomposition \cite{Oseledets2009, Oseledets2011} is a powerful tool, since Tucker decomposition suffers from the curse of dimensionality. Specifically, the Tucker core has an equal number of dimensions with the initial array, thus the complexity of the decomposition increases exponentially as $\mathcal{O}\left(dnr+r^d\right)$, with $n$ being the dimensions of the array and $r$ the Tucker ranks. In the case of TT, the decomposition of any multidimensional array consists only from matrices and 3D tensors, thus, the memory complexity remains small. However, for 3D problems, $d=3$, like the one under-study, Tucker decomposition requires less memory than the traditional TT approach. 

\renewcommand{\thefigure}{1}
\begin{figure}[ht!]
	
	\centering	
	\subfloat[]{
	\centering
	\includegraphics[scale=0.9]{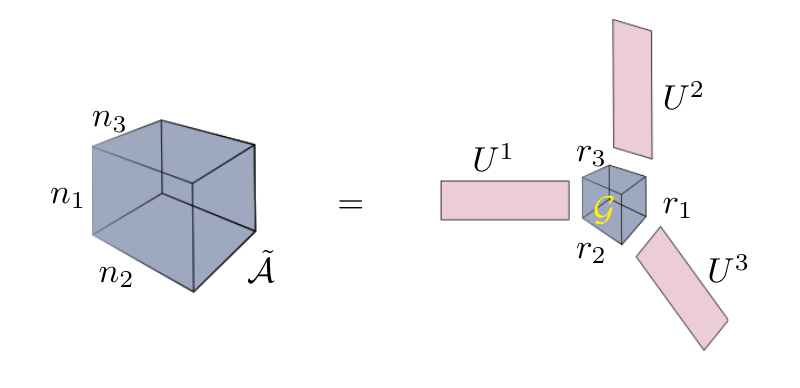}}
	
	\subfloat[]{
	\centering
	\includegraphics[scale=0.9]{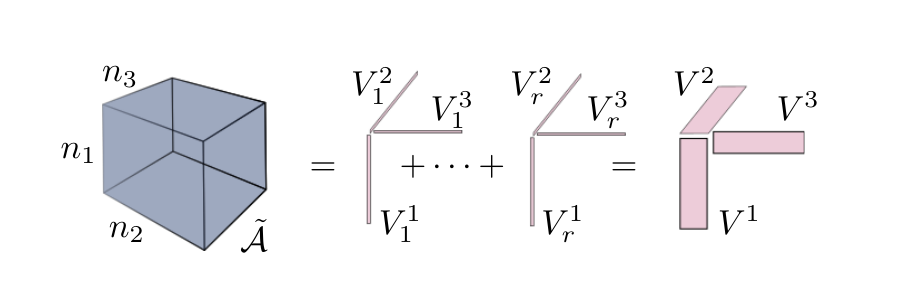}}
	
\caption{Basic tensor decompositions. (a) Tucker decomposition and (b) canonical polyadic decomposition of a low multilinear rank tensor $\Aop$.}\label{fig: n1}
\end{figure}

A different, but equally important decomposition for multidimensional arrays is the canonical polyadic model. CP decomposition is the equivalent of the SVD in two dimensions, since it approximates the original tensor with the sum of rank-one tensors.  i.e.,
\begin{equation}
\Aop \approx \tilde{\Aop} = \sum\limits_{l=1}^r V_l^1 \odot V_l^2 \odot V_l^3 
\label{eq:n15}
\end{equation}
where $V^{i} \in \mathbb{C}^{n_{i} \times r}, i =1,2,3$ are the so-called \textit{CP factors} of $\Aop$ and $r$ is the minimal number of terms to be used (\textit{canonical rank}). The outer product $\odot$ is defined in Appendix I. The CP is visualized in Fig. 1.(b). It is unique (under mild conditions \cite{Sidiropoulos2000}), however in some cases CP decomposition is ill-posed \cite{haastad1990tensor}, thus an optimal algorithm might not exist at all for a fixed number of terms \cite{Silva2008}. 

\section{Compression Algorithms} \label{sc:IV}

\subsection{Higher Order SVD}

To derive the Tucker factors $U^{1,2,3}$ and the Tucker core $\Gop$ we can use the \textit{HOSVD} algorithm proposed in \cite{Lathauwer2000}, with a proven upper error bound on the approximation. The algorithm is described briefly below.

\begin{breakablealgorithm}
\caption{HOSVD}\label{al:n1}
\begin{algorithmic}[1]

	\State $\Aop \in \mathbb{C}^{n_1 \times n_2 \times n_3}$
	
	\State Get the unfoldings matrices of $\Aop$: (Fig. 2.)
	\Statex	$A^{\left(1\right)} \in \mathbb{C}^{n_1 \times \left(n_2 \cdot n_3\right)}$, 
	\Statex $A^{\left(2\right)} \in \mathbb{C}^{n_2 \times \left(n_1 \cdot n_3\right)}$, 
	\Statex $A^{\left(3\right)} \in \mathbb{C}^{n_3 \times \left(n_1 \cdot n_2\right)}$
	
	\State Set desired tolerance $\epsilon$		
		
	\For{i=1,2,3} truncated $\text{SVD}\left(A^{\left(i\right)}\right)_{r_i}$, $\sigma_{r_i} < \frac{\epsilon}{\sqrt{3}} \sigma_{\text{max}}$ 
		\Statex $A^{\left(i\right)} \approx U^i_{r_i} \Sigma^i_{r_i} \left(V^i_{r_i}\right)^*,\: r_i = \text{rank} \left(A^{\left(i\right)} \right)$
	\EndFor	

	\State Tucker Factors: $U^i = U^i_{r_i}$

	\State Tucker Core: $\Gop = \Aop \times_1 \left(U^{1}\right)^* \times_2 \left(U^{2}\right)^* \times_3 \left(U^{3}\right)^*$ 	

\end{algorithmic}
\end{breakablealgorithm}

\noindent
The resulting tensor $\tilde{\Aop}$ satisfies the quasi-optimality condition
\begin{equation}
\norm{\Aop - \tilde{\Aop}} \leq \sqrt{3} \norm{\Aop - \Aop_{\text{best}}}
\label{eq:n16}
\end{equation}
where $\Aop_{\text{best}}$ is the best Tucker approximation of $\Aop$. Since the algorithm is based on SVD, it is stable and the low multilinear rank approximation always exists.
\par
Tucker decomposition can be constructed by using only some rows, columns and fibers of $\Aop$ with a cross-Tucker approximation algorithm \cite{Oseledets2008, Rakhuba2015} with linear complexity $\mathcal{O}\left( nr_1r_2r_3 \right)$ over the $\mathcal{O}\left( n^3\right)$ complexity of HOSVD ($n$ denotes tensor's linear size). Such algorithms are based on well-known 2D cross approximation methods \cite{Goreinov1997, Goreinov2010, Bebendorf2000, Zhao2005}.

\renewcommand{\thefigure}{2}
\begin{figure}[ht!]
\begin{center}

	\centering
	\includegraphics[scale=0.9]{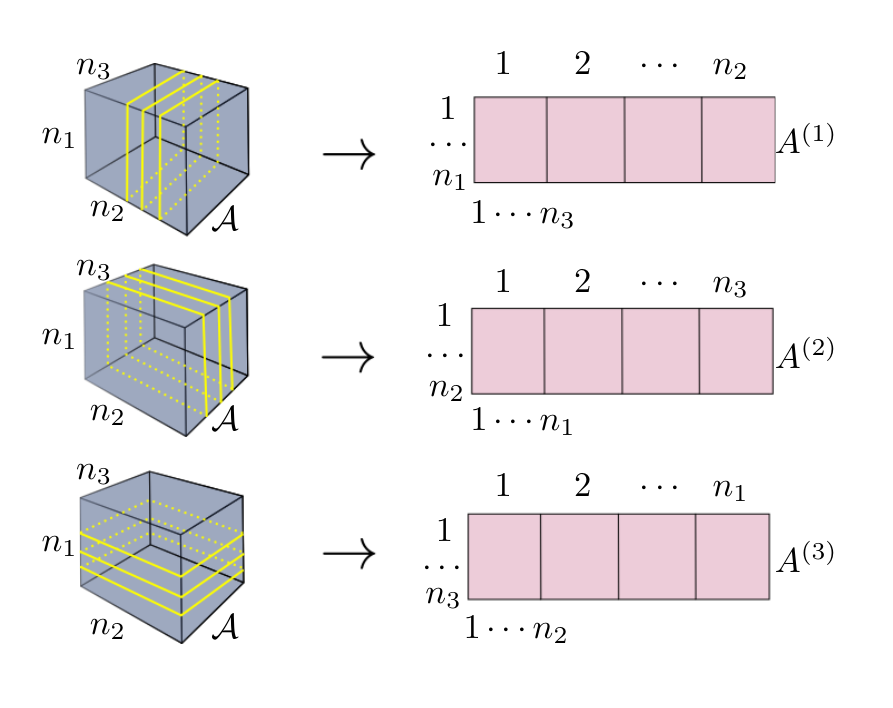}
	
\caption{Unfoldings of $\Aop$.}
\label{fig: n2}
\end{center}
\end{figure}

\subsection{Tucker + CP Method}

Another standard algorithm is the \textit{Tucker+CP method} \cite{Zhou2014,Bro1998}. In this method we use the Tucker model as a first step, which is significantly faster and more stable than CP, to compress the initial array and then, as a second step, we implement a CP algorithm on the Tucker core as shown in Fig. 3. The new decomposed form consists of 3 Tucker+CP factors. For the sake of completeness, we provide a short description of the algorithm below.

\begin{breakablealgorithm}
\caption{Tucker+CP}\label{al:n2}
\begin{algorithmic}[1]

	\State $\Aop \in \mathbb{C}^{n_1 \times n_2 \times n_3}$
	
	\State HOSVD on $\Aop \approx \Gop \times_1 U^{1} \times_2 U^{2} \times_3 U^{3}$
	
    \State Choose appropriate CP algorithm \cite{vervliet2016tensorlab} 	
	
	\State CP on $\Gop \approx \sum_{l=1}^r V_l^1 \odot V_l^2 \odot V_l^3 $		
		
	\State New Tucker+CP factors: $W^i = U^i V^i, i = 1,2,3$ 	

\end{algorithmic}
\end{breakablealgorithm}

\noindent
Typically, the step of the CP approximation is implemented by means of the alternating least squares method \cite{vervliet2016tensorlab}. Finally the approximation of $\Aop$ is given by the following equation:
\begin{equation}
\tilde{\Aop} = \sum_{l=1}^r W_l^1 \odot W_l^2 \odot W_l^3. 
\label{eq:n17}		
\end{equation}

\renewcommand{\thefigure}{3}
\begin{figure}[ht!]
\begin{center}

	\centering
	\includegraphics[scale=0.8]{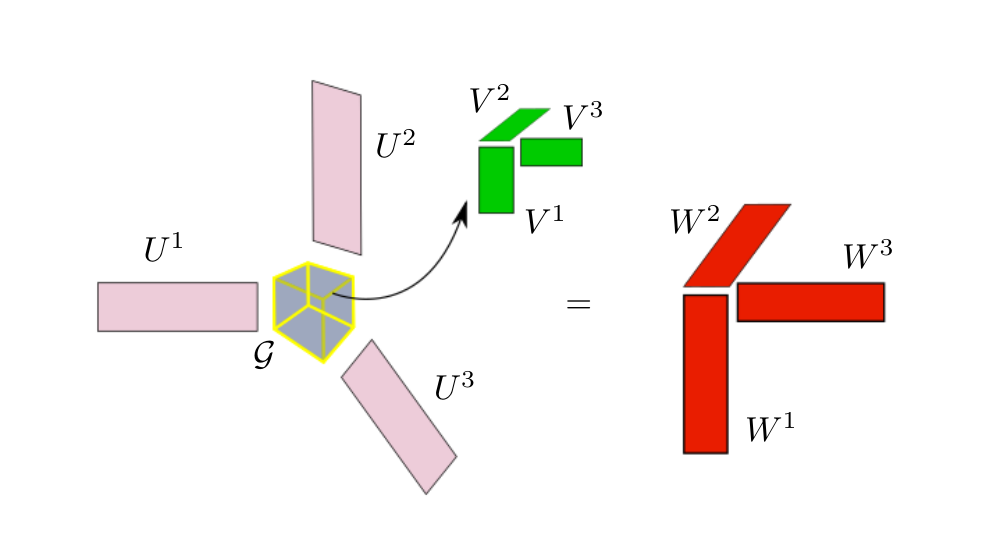}	
	
\caption{Visual representation of Tucker + CP method.}
\label{fig: n3}
\end{center}
\end{figure}

\section{Operation and Memory Footprint}  \label{sc:V}

\subsection{Circulant Embedding and Fourier Transform on the Compressed Forms}

To implement the FFT-based accelerated matrix-vector product, we need to compute the Fourier transform of each unique circulant component of $\Nop$ and $\Kop \in \mathbb{C}^{n_1 \times n_2 \times n_3}$, with $\Oop \left( N_v \log N_v \right)$ complexity. Instead, we can apply the circulant embedding and the FFT directly to the decomposed tensor forms, since these transforms do not affect their low-rank properties \cite{Rakhuba2015}. They are equivalent to univariate respective 1D transforms along the columns of the Tucker or Tucker+CP factors. The overall complexity is $\mathcal{O} \left( r_i n_i \log n_i \right)$, where $i \in \{1,2,3\}$, for each factor respectively. 

\subsection{Matrix-Vector Product Cost}

The tensor decompositions mentioned in the previous sections are used in each component of $\Nop$ and $\Kop$ in order to significantly compress them and drastically reduce the required memory for their storage. Our goal is to accelerate the matrix-vector products that arise in the iterative solution of the linear system by implementing it inside the limited memory of GPU. 
\par
Without any decomposition, the matrix-vector product requires the implementation of element-wise products between the Green's function tensors and the tensors of the unknown currents, along with multidimensional FFTs with $\Oop\left(N_v\right)$ and $\Oop\left(N_v \log N_v\right)$ complexity respectively. By using a compression algorithm, we have each operator component stored in a compressed form, thus the implementation of the product requires an extra step, that can be accomplished in two ways. We can, either decompress the appropriate component using the equations (\ref{eq:n14}) or (\ref{eq:n17}), implement the element-wise product afterwards and repeat the process for the following components by using the same buffer in memory, or, we can avoid storing any extra tensor by doing the decompression and the element-wise product simultaneously via a 6D (HOSVD) or 4D loop (Tucker+CP). The FFTs are performed as in the traditional implementation for the vector (or tensor) of the unknowns. For the HOSVD the first method is much faster, but it requires an additional buffer of memory, which for large problems, might lead to memory overflow in GPU. To avoid such impasse the second method can be used with the Tucker+CP decomposition. Not only the costly 6D loop of HOSVD is reduced to a faster 4D loop (considering that the canonical rank is small enough, $r = \min\{r_1,r_2,r_3\}$), but in addition we do not need extra memory for the buffer. We summarize both the memory footprint and the complexity of the operations of these matrix-vector product methods (the FFT cost is excluded) in Table II, for the case of FFT-JVIE, where the current is expanded with PWL functions. The Tucker and the canonical rank are $\{r_1,r_2,r_3\}$ and $r$ respectively.     

\begin{table}[!ht]
\caption{Complexity of Element-wise Product} \label{tb:n2} \centering
{\def\arraystretch{2}\tabcolsep=3pt

\begin{tabular}{ c|c|c|c }
\hline
\hline
Algorithm                  & Method                  & Operations    & Memory \\
\hline
        -                  & Traditional FFT-JVIE    & $\Oop\left( N_v \right)$          & $72N_v$ \\
\hline
\multirow{2}{*}{HOSVD}     & Component decompression & $\Oop\left( r_3N_v \right)$       & $13N_v$ \\     
						   & Multidimensional loop   & $\Oop\left( r_1r_2r_3N_v \right)$ & $12N_v$ \\
\hline
\multirow{2}{*}{Tucker+CP} & Component decompression & $\Oop\left( rN_v \right)$         & $13N_v$ \\     
						   & Multidimensional loop   & $\Oop\left( rN_v \right)$         & $12N_v$ \\
\hline
\hline
\end{tabular}
}
\end{table}        

\section{Numerical Results} \label{sc:VI}

\subsection{Multilinear Rank Study}

In this section, we first validate the linear dependence of the multilinear rank of Green's function tensors with respect to the frequency. Specifically, we present the multilinear rank of the $\Nop$ and $\Kop$ components, expanded with PWL functions, along with a simple Green's function tensor $\Gop$ (Operator's kernel $g(\vec{r}-\vec{r}')$), expanded with PWC functions. Since the multilinear rank is a triplet and we are dealing with $60 \Nop$ and $30 \Kop$ tensor components, we depict only the worst case scenario (maximum rank of all Tucker ranks of all components). The simulation is performed for a cube with unit edge length and the multilinear ranks are calculated with the HOSVD algorithm and for SVD tolerance $\epsilon = 10^{-8}$. The frequency sweep analysis is implemented for the cases $0.3i \text{ GHz }, i = 1,\dots,10$ and for three different discretizations with resolutions $\nicefrac{\lambda}{10}$, $\nicefrac{\lambda}{20}$ and $\nicefrac{\lambda}{30}$. The results are summarized in Fig. 4. The dependence between the ranks and the frequency is almost linear and the overall compression factor is remarkable.  

\renewcommand{\thefigure}{4}
\begin{figure}[ht!]
\begin{center}
\includegraphics[scale=0.59]{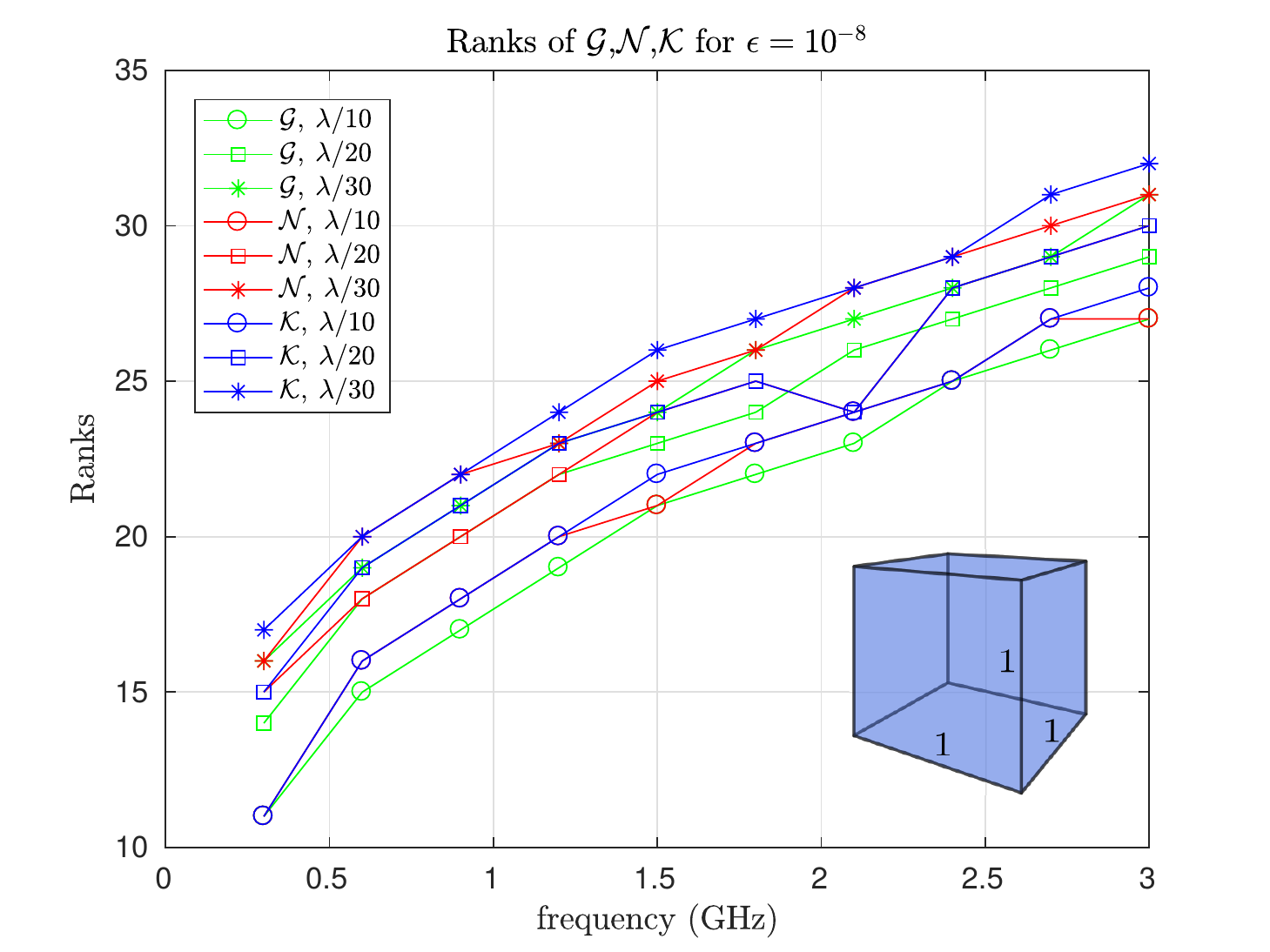}
\caption{The maximum ranks of Green's function tensors. The ranks scale almost linearly with respect to frequency.}\label{fig:n4}
\end{center}
\end{figure}   

In the following, we consider an example of a domain needed for a realistic human body model, relevant to $7$ Tesla MRI applications (the working frequency is $f \approx 298 \text{MHz}$). The domain is a rectangular cuboid with dimensions $0.538 \times 0.28 \times 1.802 \; \text{m}^3$ with a resolution of $2$ mm. The discretized version of the domain requires $270 \times 141 \times 902$ voxels and the overall memory for the unique circulant forms of $\Nop \: \& \: \Kop$ with PWL functions is $245.6$ and $122.8$ GB respectively. In Table III we provide the overall storage memory needed for the decomposed $\Nop \: \& \: \Kop$ unique components, as they are given by the HOSVD algorithm, for various tolerances. 

\begin{table}[!ht]
\caption{Memory (in MB) for the decomposed $\Nop$ and $\Kop$ \hspace{\textwidth} tensors with the hosvd} \label{tb:n3} \centering
{\def\arraystretch{2}\tabcolsep=3pt

\begin{tabular}{ c|c|c|c|c|c|c|c }
\hline
\hline
$\epsilon$ & $10^{-4}$ & $10^{-5}$ & $10^{-6}$ & $10^{-7}$ & $10^{-8}$ & $10^{-9}$ & $10^{-10}$ \\
\hline
$\Nop$     & $18.04 $ & $24.16 $ & $32.10 $ & $40.94 $ & $50.82 $ & $60.90 $ & $71.71 $ \\
$\Kop$     & $12.06 $ & $16.26 $ & $21.24 $ & $26.15 $ & $31.03 $ & $35.69 $ & $40.96 $ \\
\hline
\hline
\end{tabular}
}
\end{table} 

From Table III we can deduce that the achieved compression for all $\Nop$ and $\Kop$ components is of the orders of magnitude. In fact, from hundreds of GBs, we only need to store dozens of MBs. Consequently, memory-intensive problems, like the above, can be accurately solved, via low multilinear rank approximation techniques, even on personal laptops instead of servers equipped with immense random-access memory (RAM). In Fig. 5, the dual axis chart presents, on the right axis, the overall compression factor (number of elements of the full forms divided with the number of elements of the compressed forms) of the components, and, on the left axis, their maximum ranks, with respect to the same tolerances as before, both for HOSVD and TT-SVD \cite{Oseledets2011, Oseledets2012}. The significant compression factor (e.g., $\sim 8000\times$ for $\Nop$ and $\sim 6000 \times$ for $\Kop$ components, for $\epsilon = 1e-6$ and HOSVD) allows us to fit all the decomposed $\Nop \: \& \: \Kop$ components on commercial GPU cards. In addition, it is obvious that TT requires more memory than Tucker decomposition, so its application on a 3D dimensional problem is not as efficient. Potentially, one can compress the TT scheme more by using the Tucker decomposition on each individual tensor component of TT and reach the compression achieved with Tucker decomposition. However, such an approach will require additional operations in the decompression of each tensor component, which is crucial for the time footprint of the proposed matrix-vector products.

\renewcommand{\thefigure}{5}
\begin{figure}[ht!]
\begin{center}
\includegraphics[scale=0.58]{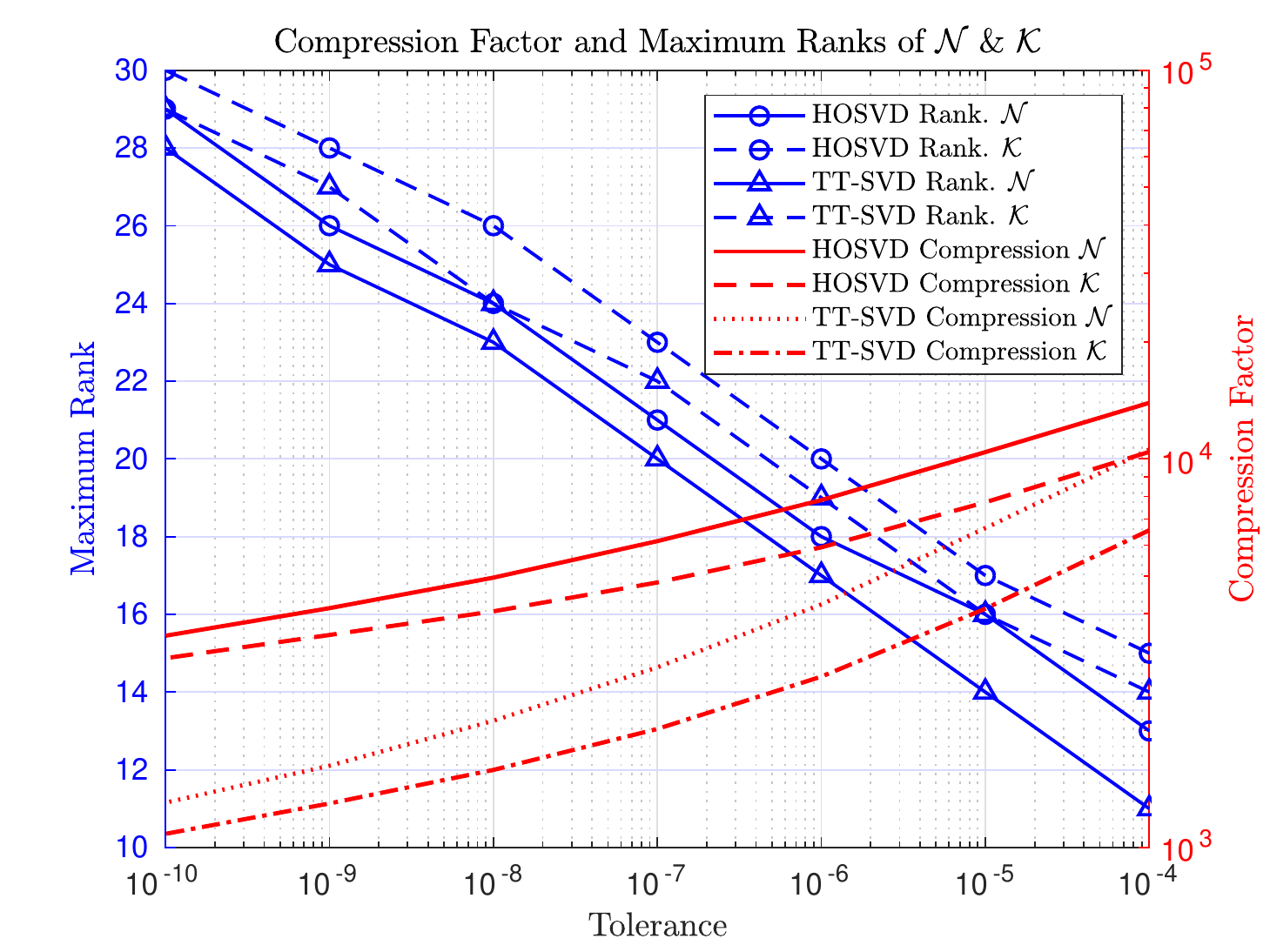}
\caption{(Left axis) Maximum rank with HOSVD and TT-SVD for the worst case scenario and (right axis) overall compression factor for all the unique components, for various tolerances, both for $\mathcal{N}$ and $\mathcal{K}$ tensors.}\label{fig:n5}
\end{center}
\end{figure}    

\subsection{Time Footprint of the Matrix-Vector Product}

Here we provide time measurements of the four novel matrix-vector product methods mentioned in Section \ref{sc:V}. We are considering double precision complex random arrays and we implement one matrix-vector product for various sizes. Specifically we use: $n = n_1 = n_2 = n_3 = 100$ or $200$ or $300$ and $r = r_1 = r_2 = r_3 = 25$, since these sizes are relevant to our simulations of interest. The calculations were done in a \textit{CentOS 6.9} operating system equipped with an \textit{Intel(R) Xeon(R) CPU E5-2699 v3 @ 2.30 GHz} with 36 cores and 2 threads per core and a \textit{NVIDIA Tesla K40M} GPU. The CPU implementation was done in C\texttt{++}, compiled with \textit{gcc 4.9.2} and the GPU one in CUDA, with \textit{nvcc V8.0.61}. The results are presented in Table IV where for completeness we include the timings measured with Matlab (version 9.2) for the same cases, both for CPU and GPU (denoted as MCPU and MGPU respectively). In the first two columns, we provide the computation times needed for decompression, along with the element-wise product. In the next two, we present the time needed for the element-wise product through the multi-dimensional loop technique. The MGPU case is not applicable for this technique, since we are unable to construct a vectorized code and cannot exploit properly the high parallelization structure of GPU. Finally, the last column provides the timings of the 3D FFT.
	
\begin{table}[!ht]
\caption{Time Footprint in \textit{ms} of Matrix Vector Product} \label{tb:n4} \centering
{\def\arraystretch{2}\tabcolsep=3pt

\begin{tabular}{ c|c|c|c|c|c|c }
\hline
\hline
\centering \multirow{2}{*}{$n,r$} & \centering \multirow{2}{*}{Case} & \centering HOSVD  & \centering Tucker+CP   & \centering HOSVD & \centering Tucker+CP  & \multirow{2}{*}{FFT} \\[0.1em]
								  &                                  & \centering decomp. & \centering decomp. & \centering loops & \centering loops &              \\[0.3em]
\hline
\multirow{4}{*}{$100,25$} & \centering C\texttt{++}   & $ 2.6  $ & $ 1.3  $ & $ 2.4e\text{+}3 $ & $ 3.3           $ & $ 1.2  $     \\[0.3em]
	         		      & \centering MCPU           & $ 5.9  $ & $ 6.4  $ & $ 2.5e\text{+}5 $ & $ 522           $ & $ 2    $     \\[0.3em]
			              & \centering CUDA           & $ 0.96 $ & $ 0.79 $ & $ 855           $ & $ 1.28          $ & $ 0.58 $     \\[0.3em]
			              & \centering MGPU           & $ 1.37 $ & $ 7.73 $ & N/A     &  N/A    & $ 0.76 $     \\[0.3em]
\hline
\multirow{4}{*}{$200,25$} & \centering C\texttt{++}   & $ 25.2 $ & $ 22.9 $ & $ 1.7e\text{+}4 $ & $ 24.6          $ & $ 20.4 $     \\[0.3em]
			              & \centering MCPU           & $ 47.3 $ & $ 56.4 $ & $ 2e\text{+}6   $ & $ 3.2e\text{+}3 $ & $ 24.6 $     \\[0.3em]
			              & \centering CUDA           & $ 6.03 $ & $ 5.79 $ & $ 5.9e\text{+}3 $ & $ 8.9           $ & $ 5.68 $     \\[0.3em]
			              & \centering MGPU           & $ 7.09 $ & $ 13.58$ & N/A     &  N/A    & $ 7.14 $     \\[0.3em]
\hline
\multirow{4}{*}{$300,25$} & \centering C\texttt{++}   & $ 74.9 $ & $ 77.5 $ & $ 5.8e\text{+}4 $ & $ 81.6          $ & $ 55.5 $     \\[0.3em]
			              & \centering MCPU           & $ 140  $ & $ 136.8$ & $ 1.5e\text{+}7 $ & $ 1e\text{+}4   $ & $ 55.6 $     \\[0.3em]
			              & \centering CUDA           & $ 17.4 $ & $ 20.1 $ & $ 1.9e\text{+}4 $ & $ 29.4          $ & $ 26.9 $     \\[0.3em]
			              & \centering MGPU           & $ 19.82$ & $ 27.83$ & N/A     &  N/A    & $ 31.43$     \\[0.3em]
\hline
\hline
\end{tabular} 

}
\end{table}

As established by the results, the GPU programming manages to calculate the matrix-vector product extremely fast. In addition, the three methods that require $\Oop\left(rn^3\right)$ complexity (HOSVD decompression, Tucker+CP decompression, Tucker+CP loops), yield equivalent results. The loop implementation of Tucker+CP doesn't require additional memory, but it is slightly slower than the decompression cases, since they can be implemented optimally with Level 3 routines from the BLAS and cuBLAS libraries for the CPU and GPU respectively. However, as expected, it is hundreds of times faster than the HOSVD loops method. Finally, we note that the highly optimized FFTw \cite{Frigo2005} and cuFFT libraries were used for the calculation of the fast Fourier transform. The use of cuFFT over FFTw significantly accelerates the convergence of the current-based VIE solver, since it is used $13$ times in every iteration of the solver of choice ($12$ FFT of the unknown current components and $1$ inverse FFT of the result of the element-wise product). Consequently, the GPU implementation of the overall matrix-vector product of the Tucker-based FFT-JVIE solver, is expected to accelerate its solution.

\subsection{EM Analysis of a Realistic Human Head}

We study the relative error between the Tucker deco\-mpo\-sition-based and the traditional FFT-JVIE implementation. Since this work is motivated by MRI applications, we illuminate a highly inhomogeneous realistic human head model with a linear polarized plane wave $\Ei = \hat{x} e^{-\iu k_0 z}$. The model is a voxelized head of the ``Duke'' body model that is provided by the Virtual Family \cite{VirtualFamily} with corresponding relative dielectric permittivity and conductivity, as illustrated in Fig. 6, for the case of a $7$ Tesla MR scanner. The domain is a cuboid, discretized with $93 \times 119 \times 125$ voxels with $2$ mm resolution and the currents are expanded with PWL functions. 

\renewcommand{\thefigure}{6}
\begin{figure}[ht!]
	
	\subfloat[]{
	\centering
	\includegraphics[scale=0.29]{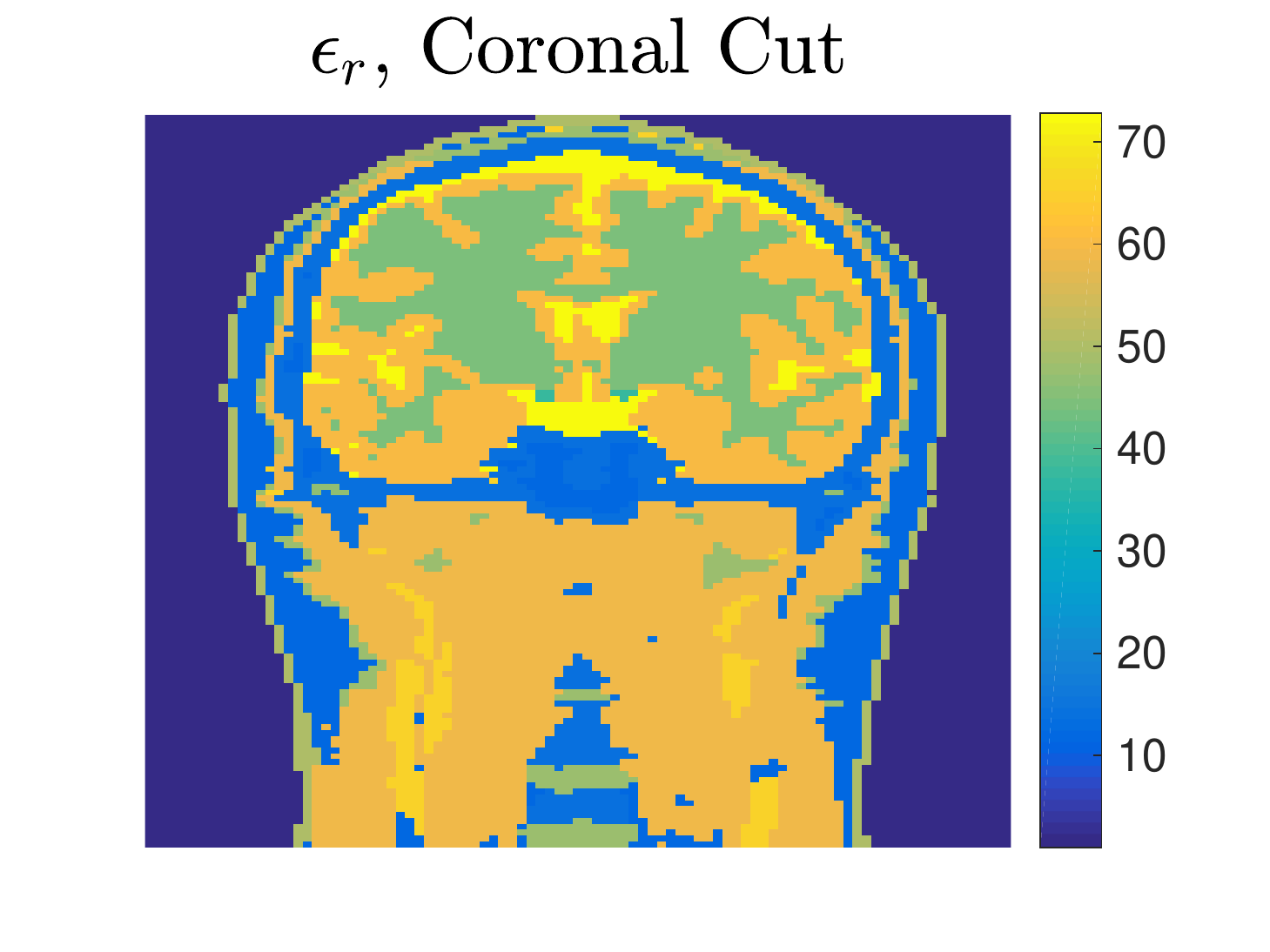}
	\includegraphics[scale=0.29]{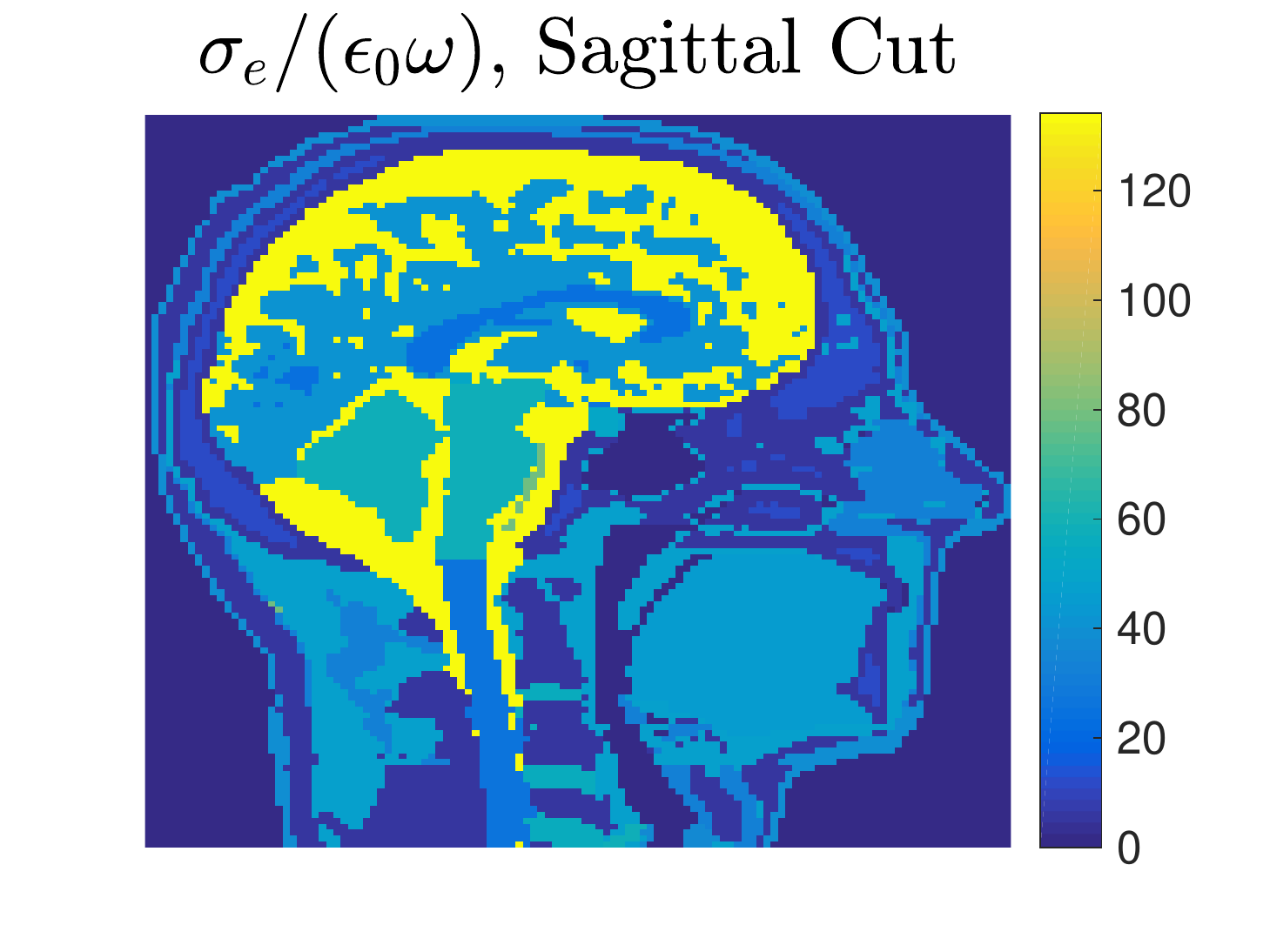}}
	
	\subfloat[]{
	\centering
	\includegraphics[scale=0.29]{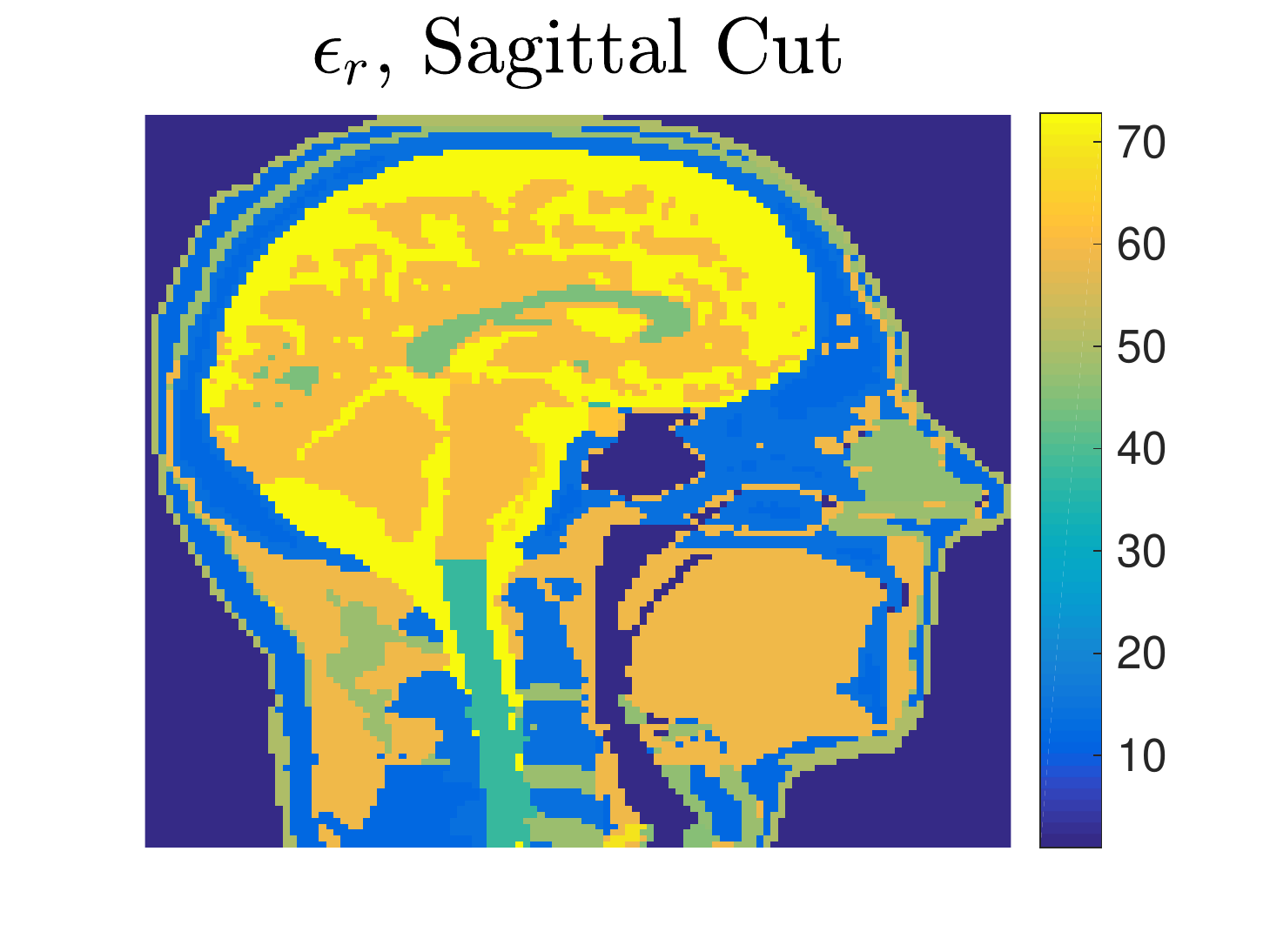}
	\includegraphics[scale=0.29]{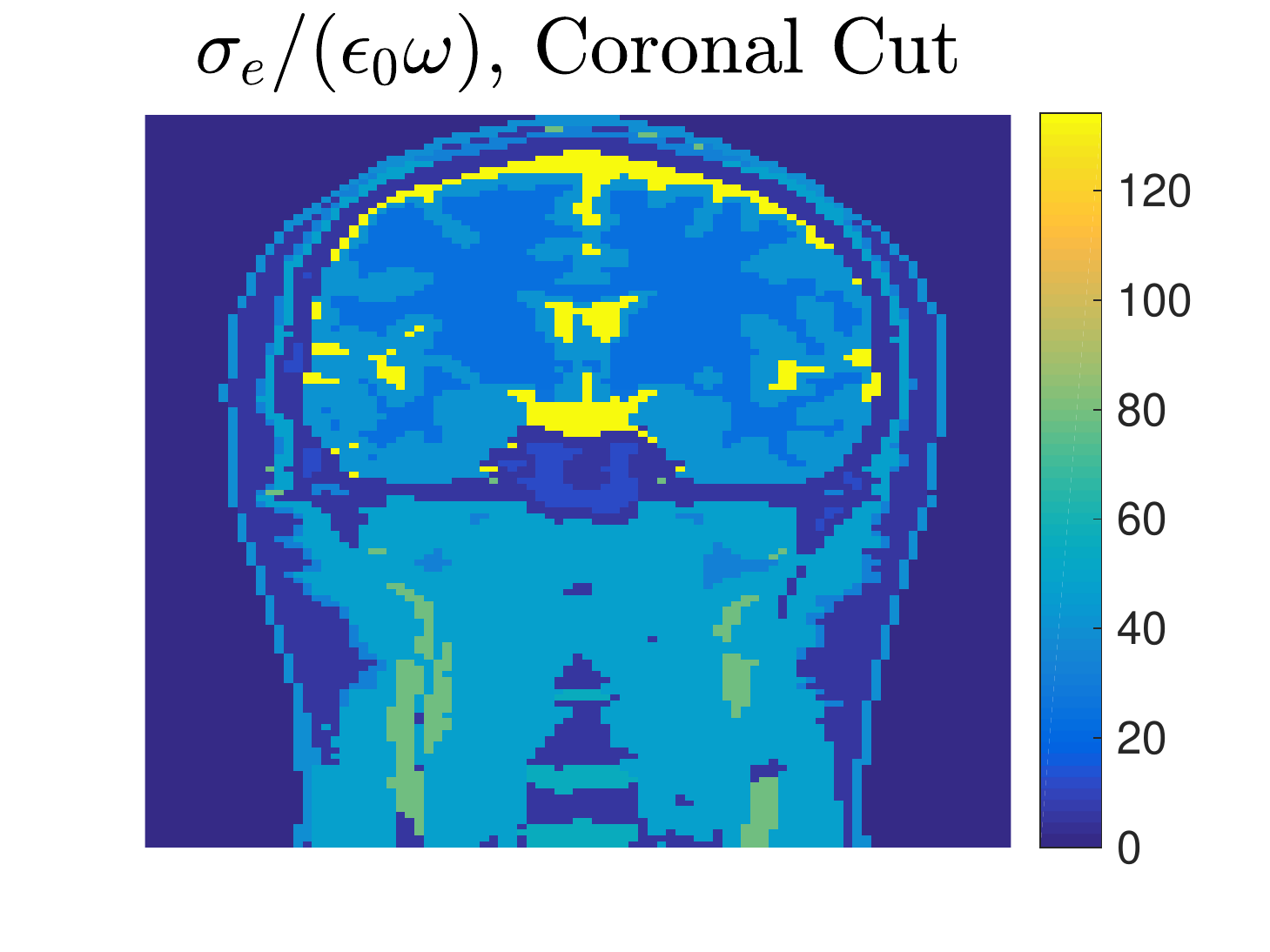}}
	
\caption{Complex dielectric permittivity of the (a) sagittal and the coronal (b) cut of the ``Duke'' realistic human head model, corresponding to $7T$ MR scanner. (From left to right) Relative permittivity $\epsilon_r$ and conductivity $\sigma_e$.}\label{fig: n5}
\end{figure} 

The calculations were performed with Matlab (version 9.2) and by embedding the aforementioned decompositions in our in-house suite MARIE \cite{MARIE}. The iterative solver of choice is the Generalized Minimum Residual method (GMRES) with tolerance $10^{-5}$ and inner and outer iterations $50$ and $200$ respectively. The traditional FFT-JVIE implementation cannot fit in the GPU, thus it is done in CPU. On the contrary, by using the HOSVD or the Tucker+CP method we are able to tackle this problem and fit our simulation in GPU and drastically reduce the convergence time of the solver. We summarize these timings in Table V.

\begin{table}[!ht]
\caption{Convergence time of the solver} \label{tb:n5} \centering
{\def\arraystretch{2}\tabcolsep=3pt

\begin{tabular}{ c|c|c }
\hline
\hline
Method      & PU  & Time in seconds \\
\hline
FFT-JVIE    & CPU & $9598$          \\
HOSVD       & GPU & $851$           \\
Tucker+CP   & GPU & $961$           \\
\hline
\hline
\end{tabular}
}
\end{table}

\par
In Fig. 7, the $L_2$ relative errors of the absorbed power $p_{\text{abs}}$ and the absolute value of the transverse magnetic flux density $\abs{b_1^+}$ between the Tucker-based decomposition implementations and the traditional FFT-JVIE are presented. These quantities hold a significant importance in MRI measurements \cite{Lattanzi2009} and are given by
\begin{equation}
p_{\text{abs}} = \frac{1}{2}\sigma_e \abs{\E}^2, \quad \abs{b_1^+} = \mu_0 \abs{\H_x + \iu \H_y}.
\label{eq:n18}
\end{equation} 
We perform our simulation for a $1000$ CP iterations and $9$ linearly scaled tolerances ($10^{-4},\dots,10^{-12}$) for the SVD.
 
\renewcommand{\thefigure}{7}
\begin{figure}[ht!]
\begin{center}
\begin{picture}(200,200)
\centering

\put(-35,0){\includegraphics[scale=0.65]{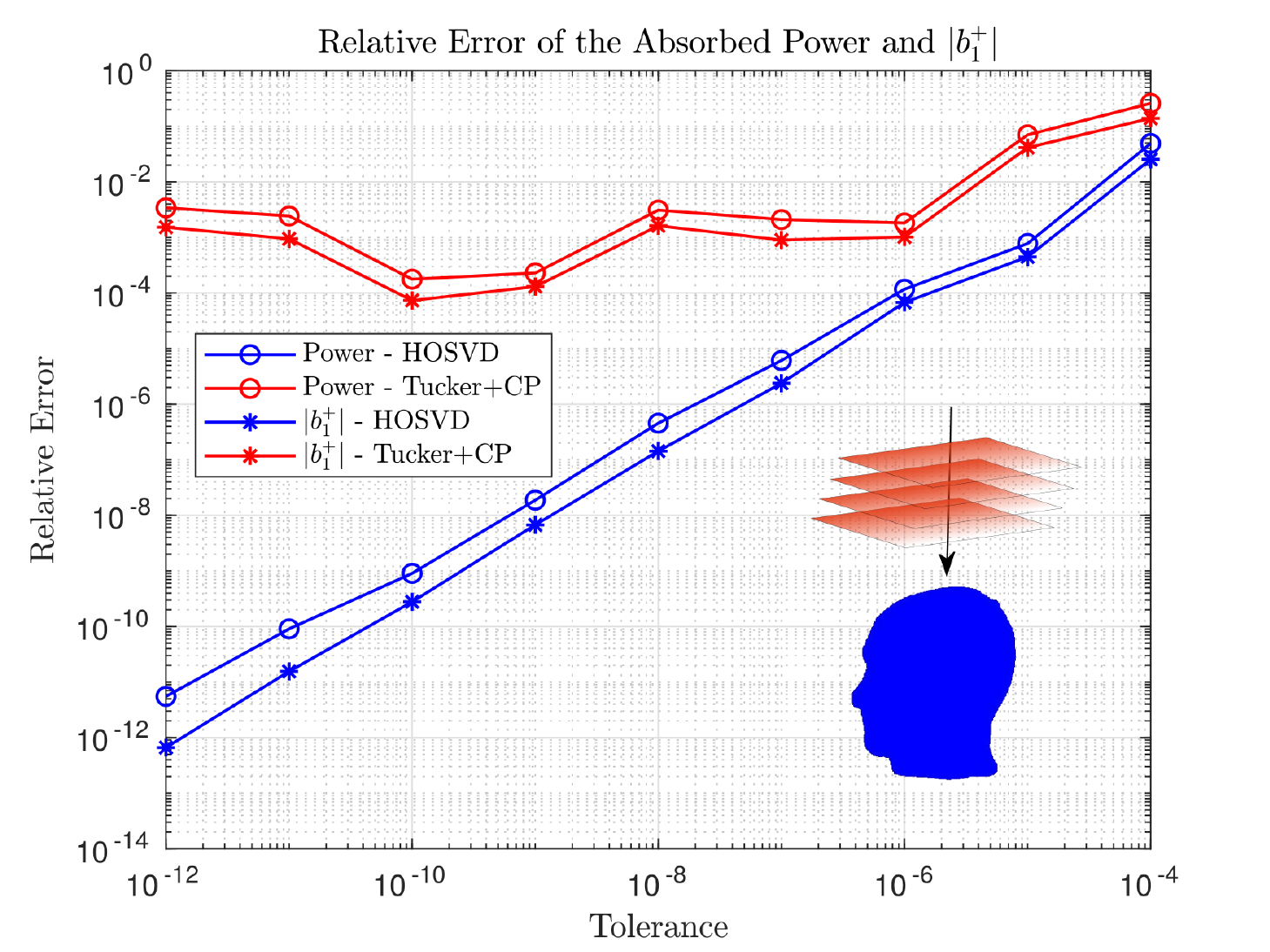}}
\put(92,75){$\Ei = \hat{x} e^{-\iu k_0 z}$}

\end{picture}
\caption{Relative error of the absorbed power and $\abs{b_1^+}$ between the Tucker-based decomposition approaches and the traditional FFT-JVIE.}\label{fig:n6}
\end{center}
\end{figure}

According to Fig. 7, the error of the HOSVD algorithm reduces linearly with respect to the given tolerance. However, Tucker+CP does not have a linear dependence with the tolerance, since the CP decomposition does not achieve the most optimal fit to the Tucker core. On the contrary, for the tolerances below $10^{-6}$, the error remains between $10^{-2}$ and $10^{-4}$ and for bigger tolerances, it is greater than $10^{-1}$. For a qualitatively validation, in Fig. 8, we illustrate the absorbed power along with $\abs{b_1^+}$ on different cuts of the head model for SVD tolerance $10^{-10}$. The values are masked outside the head for a more intuitive visualization.

\renewcommand{\thefigure}{8}
\begin{figure}[ht!]

	\subfloat[]{
	\centering
	\includegraphics[scale=0.29]{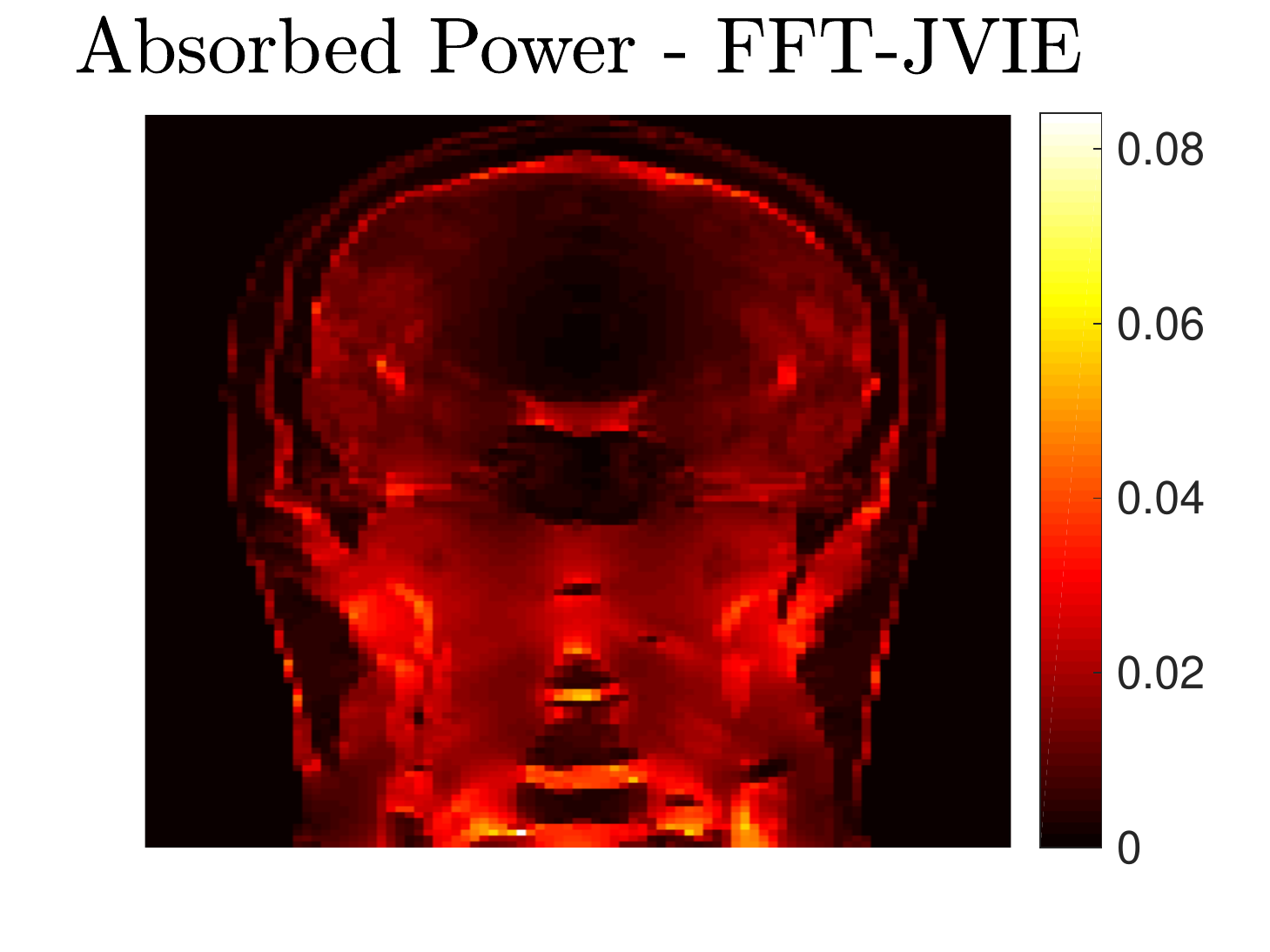}
	\includegraphics[scale=0.29]{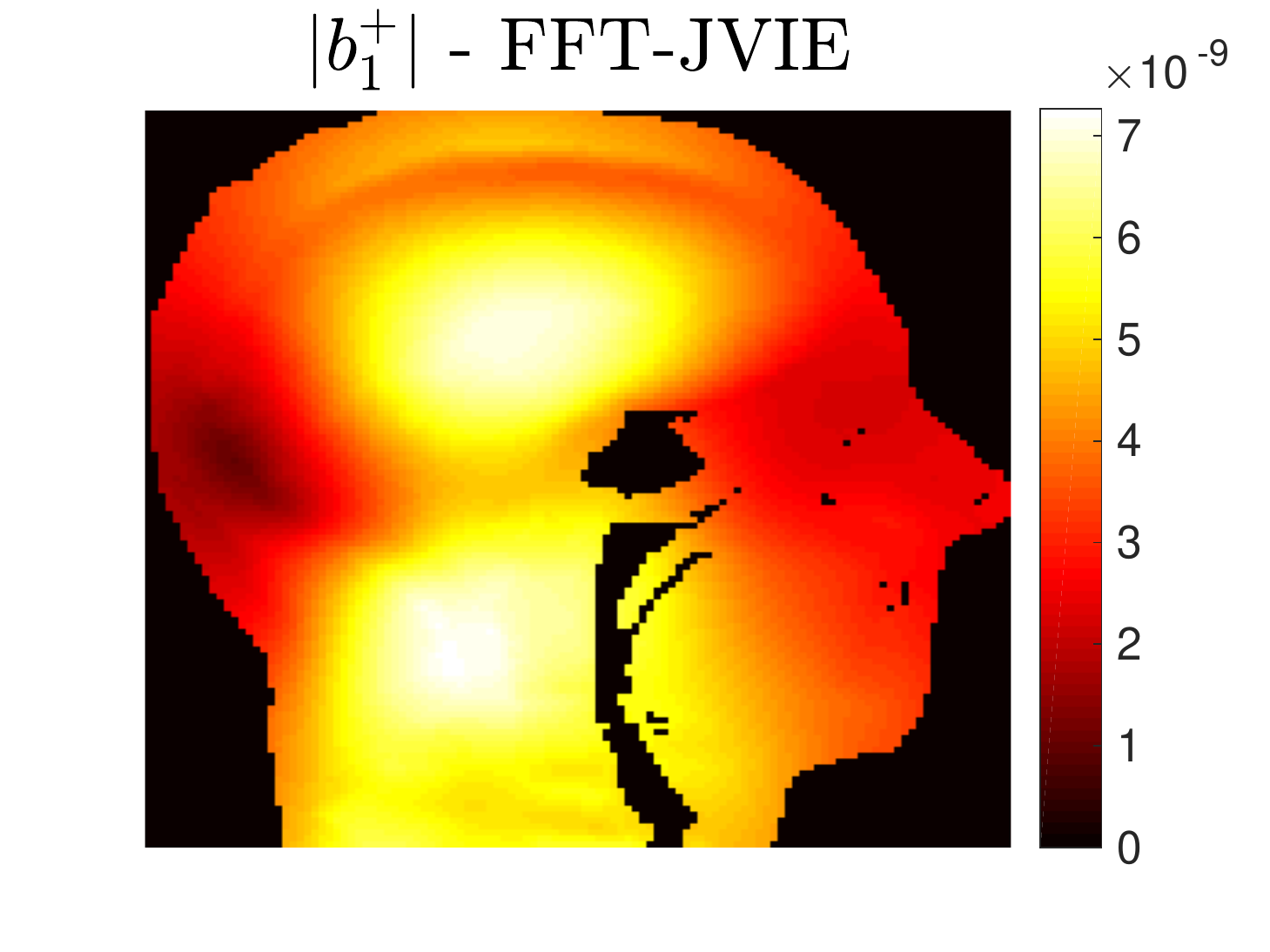}}
	
	\subfloat[]{
	\centering
	\includegraphics[scale=0.29]{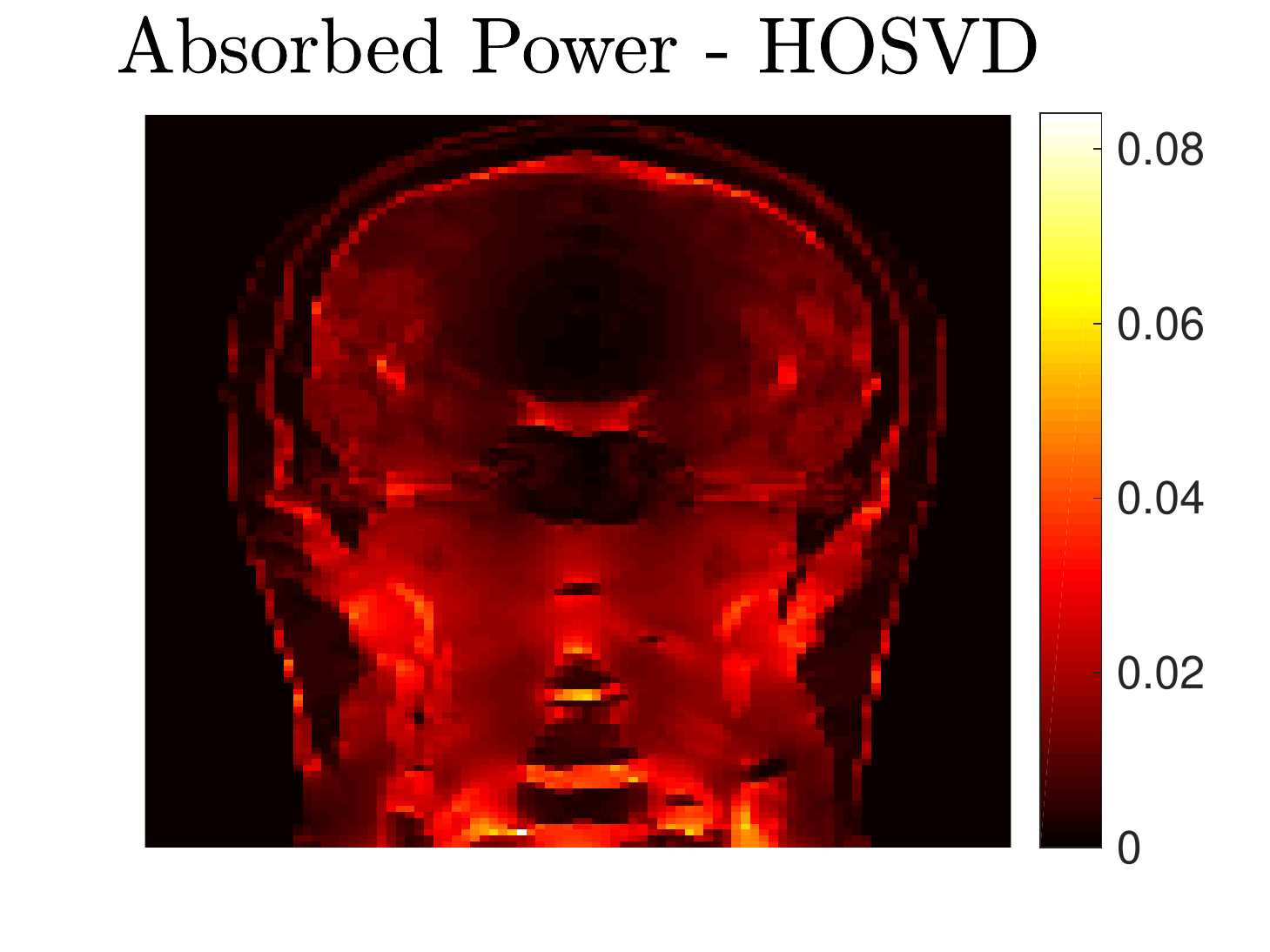}
	\includegraphics[scale=0.29]{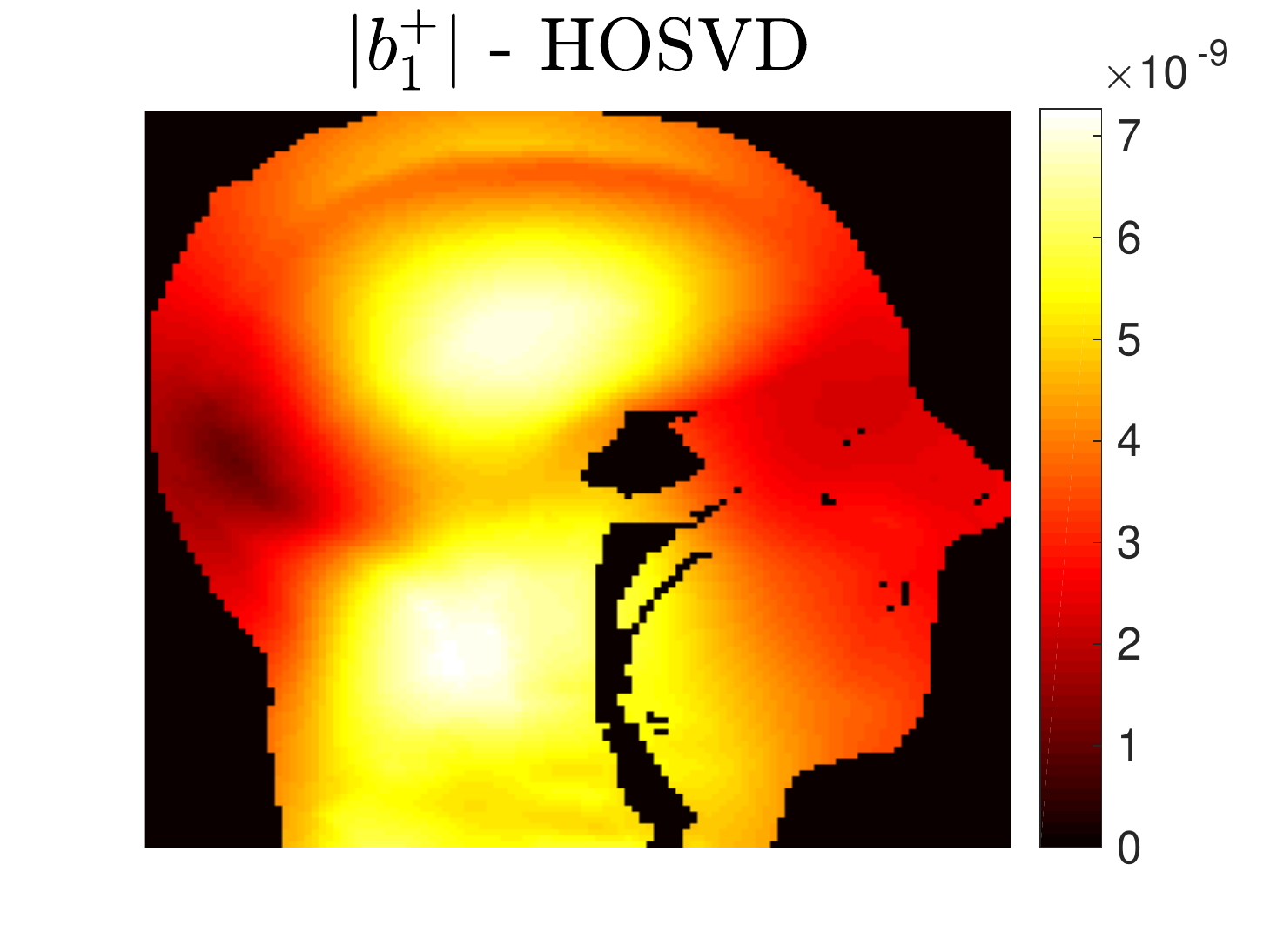}}
	
	\subfloat[]{
	\centering
	\includegraphics[scale=0.29]{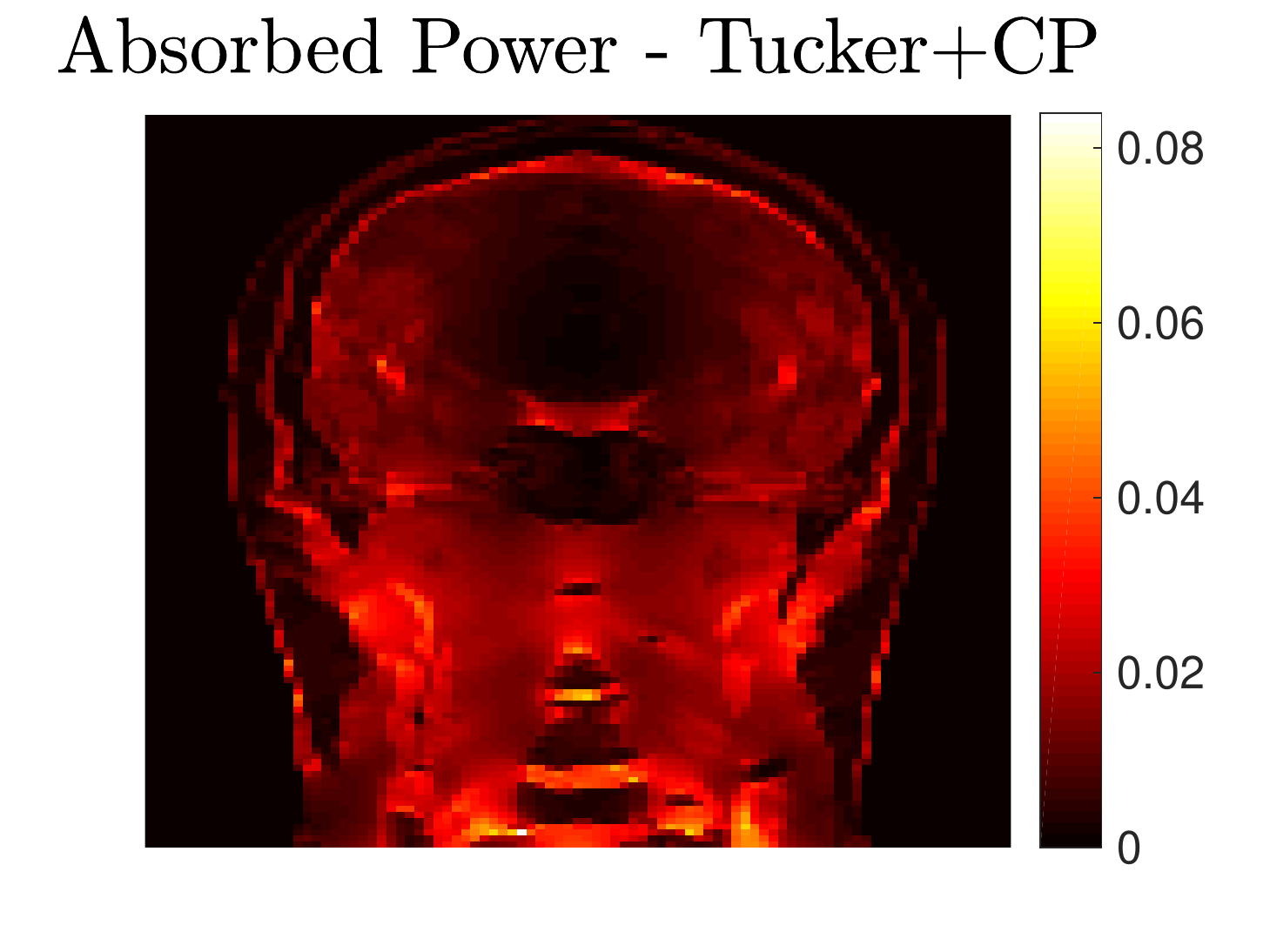}
	\includegraphics[scale=0.29]{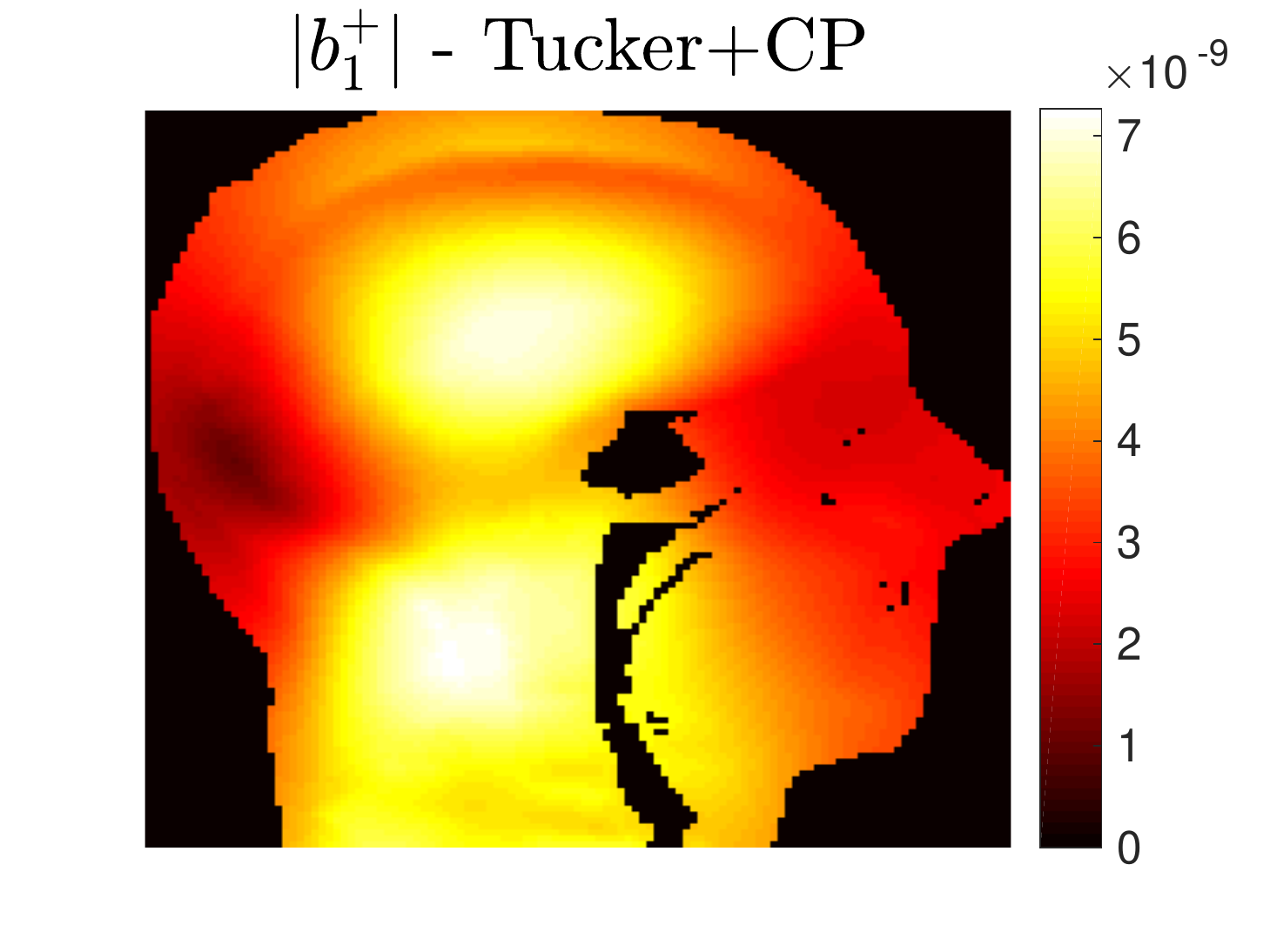}}
	
\caption{Absorbed power and $\abs{b_1^+}$ evaluated with (a) Full implementation, (b) HOSVD and (c) Tucker+CP method.}\label{fig: n6}
\end{figure}

\subsection{Scaling properties of compressed approaches}
Finally, it is of our interest to study the scaling properties of the proposed methods, namely time convergence and memory footprint as we refine the grid, thus, increasing the number of the unknowns. Concretely, a tissue-mimicking homogeneous spherical scatterer is chosen with radius $0.15$m and $\epsilon_r = 65 + \mathrm{i}\nicefrac{0.6}{\epsilon_0 \omega}$. The scatterer is  illuminated by a linear polarized plane wave $\Ei = \hat{x} e^{-\iu k_0 z}$ at $298$ MHz. The domain is a cuboid of length $0.3$m, and the voxel's isotropic resolution is $10,5,3.3,2.5,2$ mm$^3$ for each refinement of the grid, corresponding to $\sim 0.3$, $2.6$, $8.7$, $20.7$, $40.5$ million unknowns respectively, for PWL basis functions. 
\par
From Fig. 7, we can conclude that an SVD tolerance of $1e-8$ and $1000$ CP iterations are good choices, given that the GMRES tolerance is $1e-5$. Thus, for this example, we will use these settings. We perform the simulation with the traditional FFT-JVIE, the HOSVD, and the Tucker+CP approaches. In Fig. 9, on the left axis, the time convergence of the iterative solver is shown, while on the right axis we portray the memory requirements of the unique $\Nop$ components. For coarse resolutions (up to $3.3$ mm$^3$) the traditional FFT-JVIE approach can fit in the limited memory of GPU and it is the fastest method. For finer resolutions the memory needs of the unique $\Nop$ components are high, thus the simulation is forced to run in CPU, and the time footprint grows from $\sim 11$ to $\sim 486$ minutes for a $3.3$ to $2$ mm$^3$ refinement. For these resolutions, the use of compression algorithms, like HOSVD or Tucker+CP, is crucial if we want to accelerate the solution with GPU programming.  Their time footprint, e.g., for HOSVD, starts from $\sim 13$ and grows to $\sim 45$ minutes as we refine the grid (zoomed part of Fig. 9), guaranteeing an order of magnitude of acceleration.

\renewcommand{\thefigure}{9}
\begin{figure}[ht!]
\begin{center}
\includegraphics[scale=0.58]{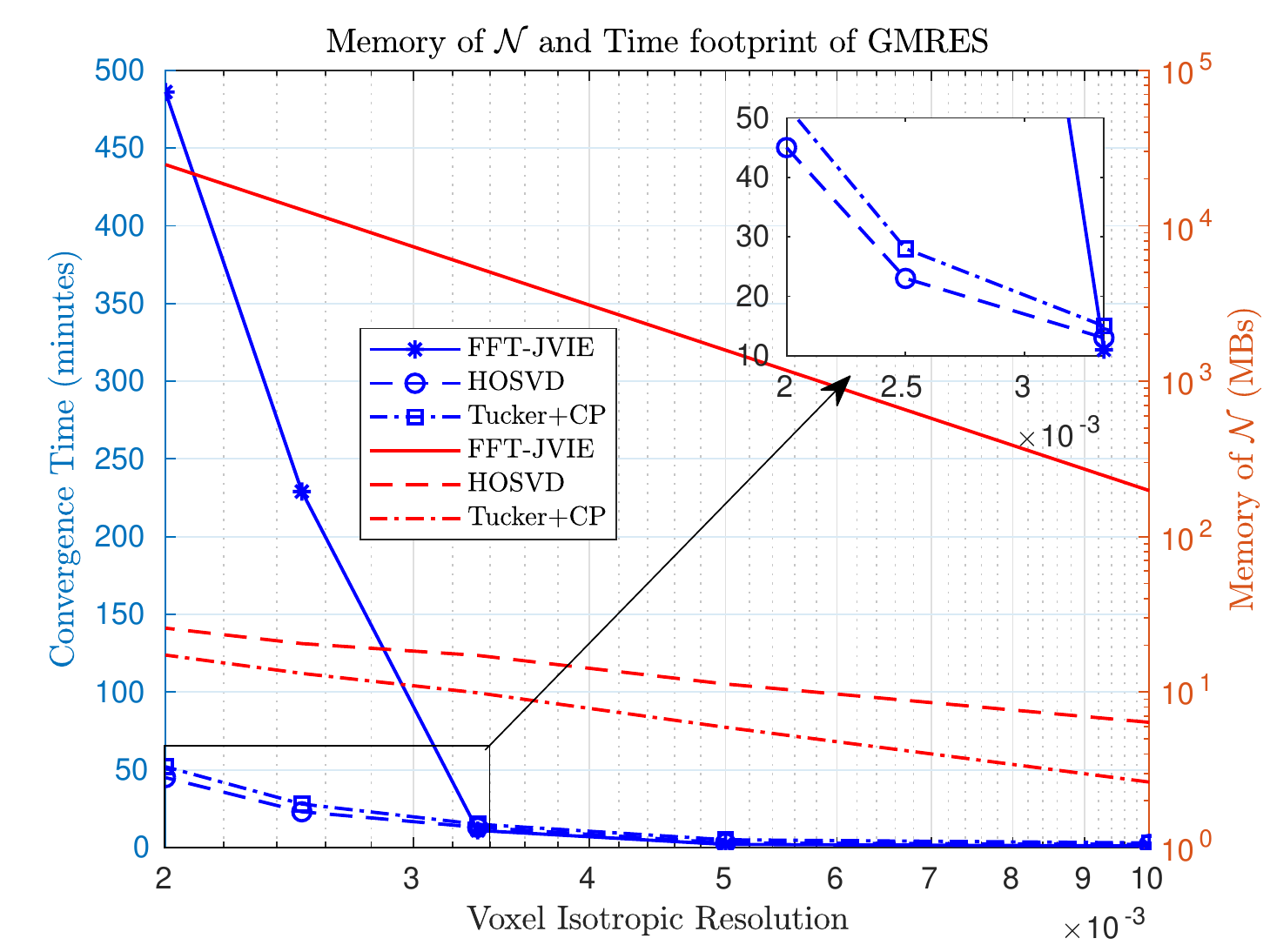}
\caption{(Left axis) Convergence time of GMRES and (right axis) memory footprint for all the unique components of $\Nop$. On the top right corner we zoom on the convergence time of finer resolutions for enhanced visualization.}\label{fig:n9}
\end{center}
\end{figure}   

In Fig. 10, we show the relative error between the scalar absorbed power, computed with the formulas in \cite{Polimeridis2015}, and the analytic solution of Mie series \cite{mie1908, MieScattering}. The error for the traditional FFT-JVIE and the HOSVD approach is identical, while the absorbed power calculated with the Tucker+CP approach does not necessarily lead to more accurate results as we refine the grid. Finally, we note that the relative error might seem high even for fine discretizations, given that we are using PWL basis functions for our simulations, although it is not surprising since the voxelized grid creates a staircase approximation of a sphere \cite{popov2002staircase}.

\renewcommand{\thefigure}{10}
\begin{figure}[ht!]
\begin{center}
\begin{picture}(200,180)
\centering
\put(-25,0){\includegraphics[scale=0.58]{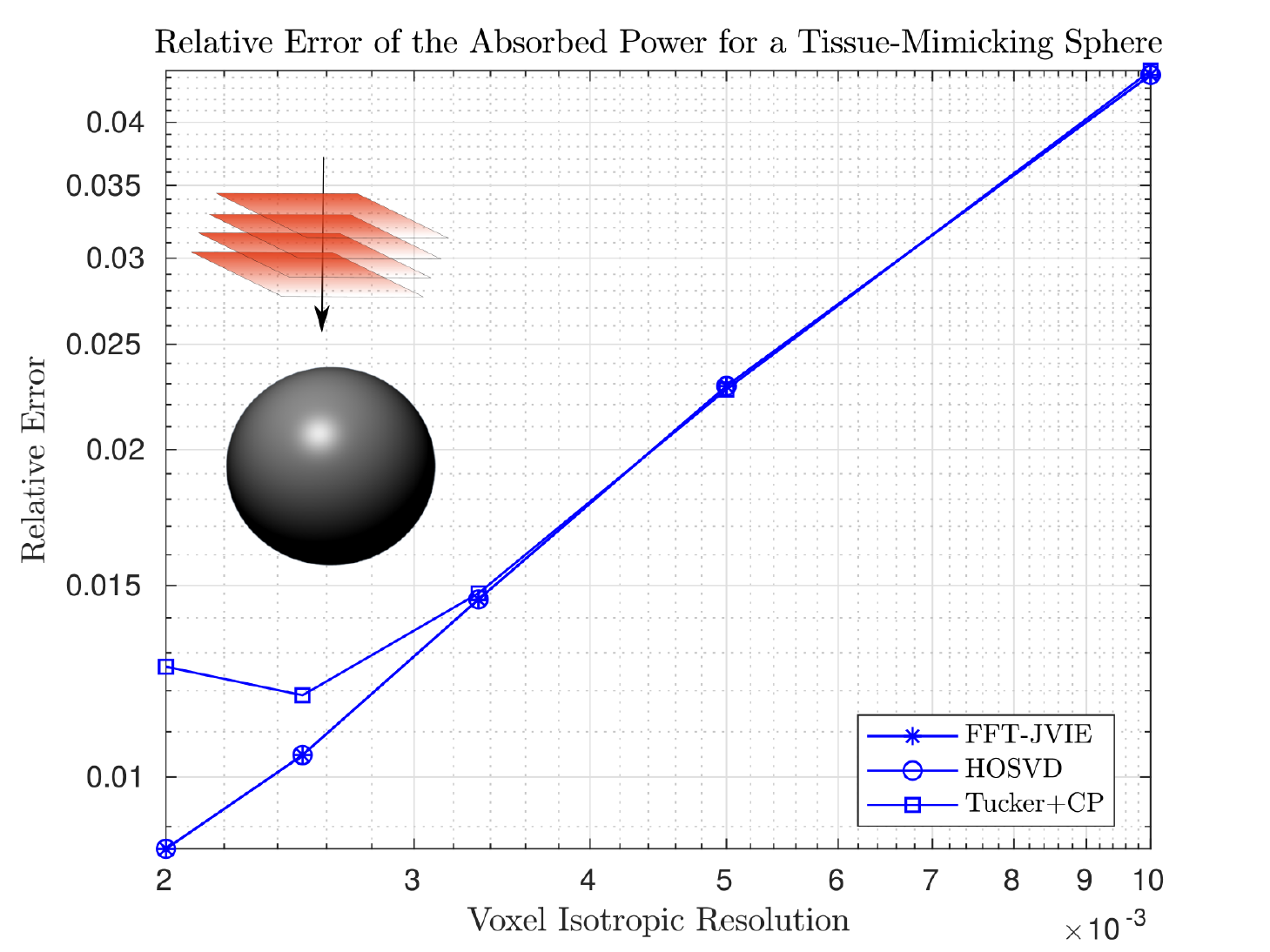}}
\put(62,135){$\Ei = \hat{x} e^{-\iu k_0 z}$}
\end{picture}
\caption{Relative error of the scalar absorbed power between the traditional FFT-JVIE, HOSVD, Tucker+CP and Mie series.}\label{fig:n10}
\end{center}
\end{figure}   

\section{Conclusion} \label{sc:VII}

In this work, a method for significant memory reduction of Green's function tensors arising in current-based VIE, via Tucker-based decompositions is presented. This allows us to handle efficiently higher order polynomial basis functions and use the highly parallel performance of GPU programming in order to accelerate significantly the numerical evaluation. Similar compression is expected for other FFT-based VIE formulations (flux and field based) since they consist of similar Green's function-based kernels. Finally, the presented work can be used to speed-up the time-consuming inverse electromagnetic scattering problems e.g., the recently proposed Global Maxwell Tomography method \cite{serralles2019noninvasive}, where the forward problem needs to be solved hundreds of times in order to retrieve an accurate dielectric property mapping of biological tissue. 

\section*{Acknowledgments}
The authors are grateful to S. P. Groth and J. K. White for valuable discussions. This work was supported by grants from the Skoltech-MIT Next Generation Program.

\renewcommand\appendixname{\textsc{\normalsize Appendix I}}
\appendix[]

\textbf{Definition A.1.} The $L_2$ inner product $\langle \cdot \rangle_V$ between two vectors $\matvec{u},\matvec{v} \in \mathbb{C}^n$ is defined as follows:
\begin{equation}
\langle \matvec{u},\matvec{v} \rangle_V = \int\limits_V \matvec{u}^*\matvec{v} dV.
\label{eq:ap1}
\end{equation}

\textbf{Definition A.2.} The Frobenius norm of a tensor $\Aop \in \mathbb{C}^{n_1 \times n_2 \times n_3}$ reads
\begin{equation}
\norm{\Aop}_F \triangleq \left( \sum_{i=1}^{n_1} \sum_{j=1}^{n_2} \sum_{k=1}^{n_3} \abs{\Aop_{i j k}}^2 \right)^{\frac{1}{2}}.
\label{eq:ap2}
\end{equation}

\textbf{Definition A.3.} The n-mode product of a tensor $\Aop \in \mathbb{C}^{m_1 \times m_2 \times m_3}$ with a matrix $U^{n} \in \mathbb{C}^{q_n \times m_n},n=1,2,3$ results to a new tensor $\Bop^{n}$ obtained by the convolution over the $n$ axis, i.e., for $\Bop^{1} \in \mathbb{C}^{q_1 \times m_2 \times m_3}$:
\begin{equation}
\Bop^{1} = \Aop \times_1 U^{1}, \quad \Bop_{i j k} = \sum_{t = 1}^{m_1} \Aop_{t j k} U^{1}_{i t}. \\
\label{eq:ap3}
\end{equation}

\textbf{Definition A.4.} The outer product $\odot$ between two multidimensional arrays $A \in \mathbb{C}^{n_1 \times \cdots \times n_N},B \in \mathbb{C}^{m_1 \times \cdots \times m_M}$ is defined as follows
\begin{equation}
\left( \matvec{A} \odot \matvec{B} \right)_{i_1,\cdots,i_N,j_1,\cdots,j_M} = A_{i_1,\cdots,i_N} \cdot B_{j_1,\cdots,j_M}.
\label{eq:ap4}
\end{equation}

\Urlmuskip=0mu plus 1mu\relax
\bibliographystyle{IEEEtran}
\bibliography{IEEEabrv,References}

\vskip -2\baselineskip plus -1fil

\begin{IEEEbiography}
[{\includegraphics[width=1in,height=1.25in,clip]{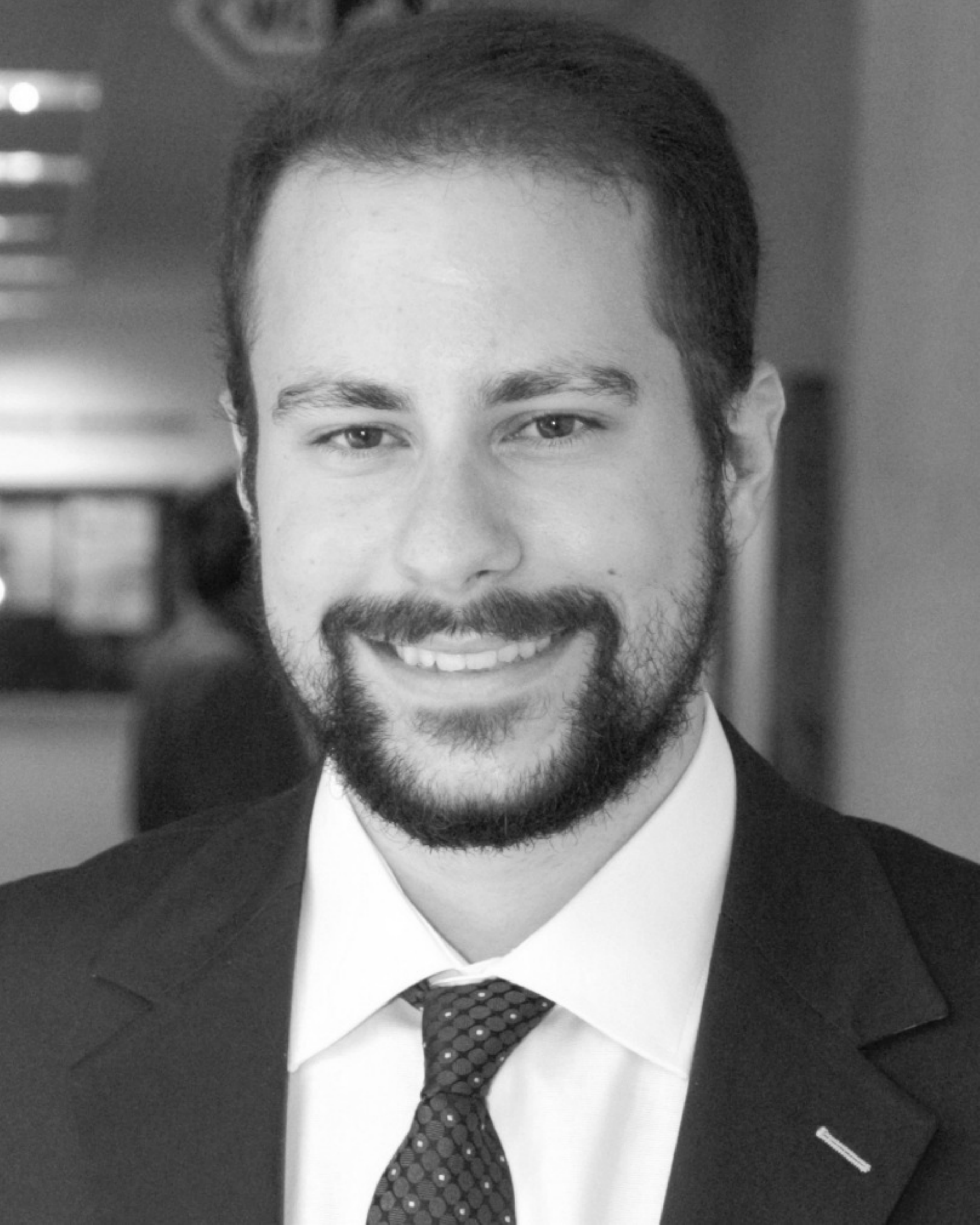}}]
{Ilias I.~Giannakopoulos} received the Diploma degree in electrical engineering and computer science from the Aristotle University of Thessaloniki, Greece, in 2016. He is currently working towards his Ph.D. degree in Skolkovo Institute of Science and Technology in the Center for Computational Data-Intensive Science and Engineering. 
\par
His research is focused on computational electromagnetics, with an emphasis on volume and surface integral equation methods, numerical linear algebra and magnetic resonance imaging.
\par
Mr. Giannakopoulos was a recipient of an honorary scholarship from the Greek State Scholarships Foundation (2014) for his excellence in his undergraduate studies for the academic year 2011-2012, and the IEEE AP-S Student Paper Competition Honorable Mention Award of the 2018 IEEE International Symposium on Antennas and Propagation.
\end{IEEEbiography}

\vskip -2\baselineskip plus -1fil

\begin{IEEEbiography}
[{\includegraphics[width=1in,height=1.25in,clip]{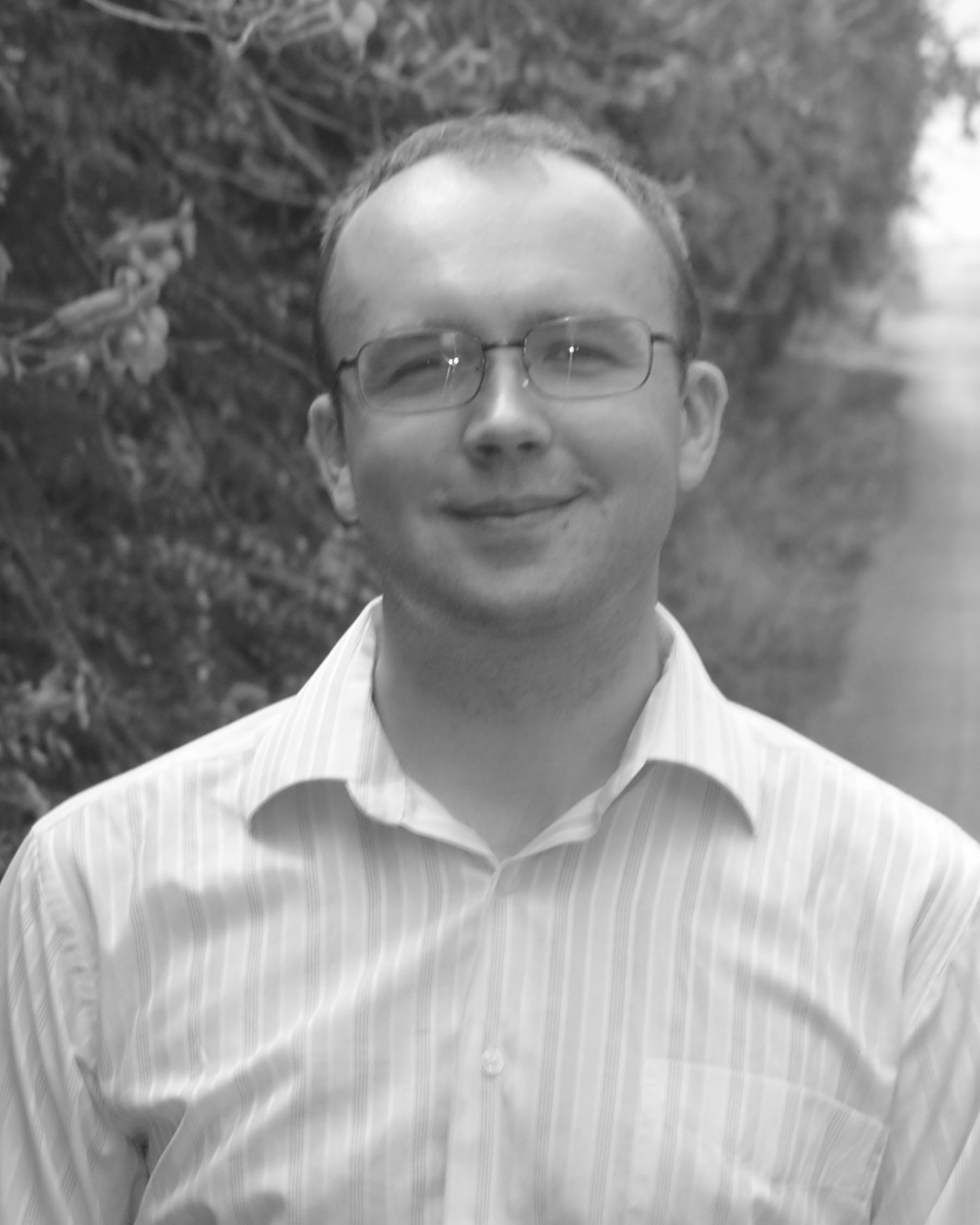}}]
{Mikhail S.~Litsarev} was born in Moscow Region, Russia, in 1983. He received the Diploma in condensed matter physics from Moscow Engineering Physics Institute in 2006 and the Ph.D. degree in theoretical physics from P.N. Lebedev Physical Institute of the Russian Academy of Sciences in 2010. 
\par
From 2011 to 2013, he was a Postdoctoral Research Associate at the Department of Physics and Astronomy, Uppsala University, Sweden, where he was a member of computational physics and condensed matter theory Group at the Materials Theory Division. Since 2013, he is a Postdoctoral Research Associate at the Skolkovo Institute of Science and Technology, Moscow, Russia. His research interests revolve around computational physics, numerical linear algebra and C++.
\end{IEEEbiography}

\vskip -2\baselineskip plus -1fil

\begin{IEEEbiographynophoto}
{Athanasios G.~Polimeridis} (SM'16) was born in Thessaloniki, Greece, in 1980. He received the Diploma and Ph.D. degrees in electrical engineering and computer science from the Aristotle University of Thessaloniki in 2003 and 2008, respectively.
\par
From 2008 to 2012, he was a Post-Doctoral Research Associate with the Laboratory of Electromagnetics and Acoustics, {\'E}cole Polytechnique F{\'e}d{\'e}rale de Lausanne, Switzerland. From 2012 to 2015, he was a Post-Doctoral Research Associate with the Massachusetts Institute of Technology, Cambridge, MA, USA, where he was a member of the Computational Prototyping Group, Research Laboratory of Electronics. From 2015 to 2018 he was an Assistant Professor with the Skolkovo Institute of Science and Technology, Moscow, Russia. He is currently the VP of Scientific Computing with Q Bio, CA, USA. His research interests include computational methods for problems in physics and engineering (classical electromagnetics, quantum and thermal electromagnetic interactions, and magnetic resonance imaging) with an emphasis on the development and implementation of integral-equation-based algorithms. 
\end{IEEEbiographynophoto}

\end{document}